\documentclass[11pt,reqno,a4paper]{amsart}
\usepackage{amssymb,amsmath}
\usepackage[latin1]{inputenc}  %FANTASTICO PER ACCENTI IN ITALIANO
\usepackage{color}

\def\enne{\mathbb{N}}
\def\R{\mathbb{R}}
\def\erre{\mathbb{R}}
\def\m1{{I\!\!M}}
\newcommand{\rdue}{\erre^2}

\def\sideremark#1{\ifvmode\leavevmode\fi\vadjust{\vbox to0pt{\vss% the remark
 \hbox to 0pt{\hskip\hsize\hskip1em%                          will appear only
 \vbox{\hsize2.1cm\tiny\raggedright\pretolerance10000%          on the side
  \noindent #1\hfill}\hss}\vbox to15pt{\vfil}\vss}}}%

\newcommand{\de}{\omega}
\newcommand{\grad}{\nabla}

\renewcommand{\to}{\rightarrow}
\newcommand{\pa}{\partial}

\newcommand{\ino}{\int_{\Omega}}
\newcommand{\E}{\exists}

\newcommand{\fo}{\forall\,}

\newcommand{\be}{\beta}

\DeclareMathOperator{\dist}{dist}

\newcommand{\supu}{\overline{u}}
\newcommand{\subu}{\underline{u}}

\newcommand{\rife}[1]{(\ref{#1})}
\newcommand{\ov}[1]{\overline{#1}}
\newcommand{\un}[1]{\underline{#1}}
\newcommand{\scp}{\scriptscriptstyle}
\newcommand{\dsp}{\displaystyle}

\renewcommand{\dfrac}{\displaystyle\frac}
\newcommand{\fineproof}{\hspace{\fill}$\square$}

\renewcommand{\i}{\infty}

\newcommand{\tsp}{\textstyle}
\renewcommand{\a}{\rho}
\newcommand{\ai}{a}
\newcommand{\bi}{b}
\newcommand{\s}{s}

\newcommand{\eps}{\varepsilon}
\newcommand{\dt}{\delta}

\newcommand{\bt}{\beta}

\newcommand{\sg}{\sigma}
\newcommand{\ga}{\gamma}
\newcommand{\om}{\Omega}
\newcommand{\lm}{\lambda}
\newcommand{\vp}{\varphi}

\newcommand{\ul}{u^{\scp(\lm)}}
\newcommand{\ull}{u_{\scp{\lm}}}
\newcommand{\vl}{v^{\scp(\lm)}}

\newcommand{\e}[1]{{\,\dsp e^{\dsp #1}}}

\newtheorem{theorem}{Theorem}[section]
\newtheorem{proposition}[theorem]{Proposition}
\newtheorem{lemma}[theorem]{Lemma}
\newtheorem{corollary}[theorem]{Corollary}
\newtheorem{remark}[theorem]{Remark}
\newtheorem{definition}[theorem]{Definition}
\newcommand{\brm}{\begin{remark}}
\newcommand{\erm}{\end{remark}}
\newcommand{\bdf}{\begin{definition}}
\newcommand{\edf}{\end{definition}}
\newcommand{\bte}{\begin{theorem}}
\newcommand{\ete}{\end{theorem}}
\newcommand{\bpr}{\begin{proposition}}
\newcommand{\epr}{\end{proposition}}
\newcommand{\ble}{\begin{lemma}}
\newcommand{\ele}{\end{lemma}}
\newcommand{\bco}{\begin{corollary}}
\newcommand{\eco}{\end{corollary}}
\newcommand{\beq}{\begin{equation}}
\newcommand{\eeq}{\end{equation}}
\newcommand{\bdm}{\begin{displaymath}}
\newcommand{\edm}{\end{displaymath}}

\newcommand{\graf}[1]{\left\{\begin{array}{ll}#1\end{array}\right.}
\newcommand{\virg}[1]{\textquotedblleft#1\textquotedblright}

\newcommand{\bel}{\begin{equation}}
\newcommand{\eel}{\end{equation}}

\begin{document}
\numberwithin{equation}{section}
\parindent=0pt
\hfuzz=2pt
\frenchspacing

\title[Supercritical Mean field equations]{Supercritical Mean Field Equations on convex domains and the Onsager's
statistical description of two-dimensional turbulence}

\author[D. Bartolucci \& F. De Marchis]{Daniele Bartolucci$^{(1,\ddag)}$\& Francesca De Marchis$^{(2,\ddag)}$}


\thanks{$^{(1)}$Daniele Bartolucci, Department of Mathematics, University
of Rome {\it "Tor Vergata"}, \\  Via della ricerca scientifica n.1, 00133 Roma,
Italy. e-mail:bartoluc@mat.uniroma2.it}

\thanks{$^{(2)}$Francesca De Marchis, Department of Mathematics, University
of Rome {\it "Tor Vergata"}, \\  Via della ricerca scientifica n.1, 00133 Roma,
Italy. e-mail:demarchi@mat.uniroma2.it}

\thanks{$^{(\ddag)}$Research partially supported by FIRB project {\sl
Analysis and Beyond} and by MIUR project {\sl Metodi variazionali e PDE non lineari}}

\begin{abstract}
We are motivated by the study of the Microcanonical Variational Principle within the Onsager's
description of two-dimensional turbulence in the range of energies where the equivalence of statistical ensembles fails.
We obtain sufficient conditions for the existence and multiplicity of solutions for the corresponding Mean Field
Equation on convex and "thin" enough domains in the supercritical (with respect to the Moser-Trudinger inequality) regime.
This is a brand new achievement since existence results in the supercritical region were previously known
\un{only} on multiply connected domains.
Then we study the structure of these solutions by the analysis of their linearized problems
and also obtain a new uniqueness result for solutions of the Mean Field Equation on thin domains whose
energy is uniformly bounded from above. Finally we evaluate the asymptotic expansion of those solutions with respect
to the thinning parameter
and use it together with all the results obtained so far to solve the Microcanonical Variational Principle in a small
range of supercritical energies where the entropy is eventually shown to be concave.
\end{abstract}
\maketitle

\tableofcontents

\section{Introduction}
\setcounter{equation}{0}

In a pioneering paper \cite{On} L. Onsager
proposed a statistical theory of two-dimensional turbulence based on the N-vortex model
\cite{New}.
We refer to \cite{ESr} for an historical review and to \cite{MP} and the introduction in \cite{ESp} for a
detailed discussion about this theory and its range of applicability in real world models.
More recently those physical arguments was turned into rigorous proofs \cite{clmp1}, \cite{clmp2}, \cite{K}, \cite{KL}.
Together with other well known physical \cite{bav}, \cite{suzC}, \cite{sy2}, \cite{T}, \cite{tar}, \cite{yang}
and geometrical \cite{cygc}, \cite{KW}, \cite{Troy} applications,
these new results were the motivation for the lot of efforts in the understanding
of the resulting mean field \cite{clmp1}, \cite{clmp2} Liouville-type \cite{Lio} equations.
We refer the reader to
\cite{B2}, \cite{bl}, \cite{BM2}, \cite{bt}, \cite{bls}, \cite{bm}, \cite{CCL},  \cite{ChK},
\cite{CL1}, \cite{cli1}, \cite{CLin3},
\cite{CLin1}, \cite{CLin2}, \cite{CLin4},  \cite{CSW}, \cite{Kwan1}, \cite{DJLW},
\cite{dj}, \cite{EGP}, \cite{KMdP}, \cite{yy},
\cite{Lin1}, \cite{lin7}, \cite{linwang},
\cite{Mal1}, \cite{Mal2}, \cite{dem}, \cite{dem2}, \cite{NT}, \cite{OS}, \cite{pt}, \cite{st},
\cite{suz}, \cite{suzB}, \cite{T3}, \cite{w}, and more recently
\cite{barjga}, \cite{BLin2}, \cite{BLin3}, \cite{BMal}, \cite{BDeM}, \cite{malru} and the
references quoted therein.\\

In spite of these efforts it seems that there are some basic questions arising in \cite{clmp2} which have
been left unanswered so far. These are our main motivations and this is why we
will begin our discussion with a short review of some of the results obtained in \cite{clmp2}
as completed in \cite{CCL}.

\bdf\label{defsimp}
Let $\om\subset\R^2$ be any open, bounded and simply connected domain.
We say that $\om$ is simple if $\pa\om$ is the support of a simple and rectifiable Jordan curve.\\
Let $\om$ be a simple domain. We say that it is regular if (see also \cite{CCL}):\\
(-) its boundary $\pa \om$ is the support of a continuous and piecewise $C^2$ curve
$\pa\om=\mbox{supp}(\gamma)$ with bounded first derivative $\|\gamma^{'}\|_{\infty}\leq C$ and
at most a finite number of corner-type points $\{p_1,\ldots,p_m\}$, that is, the inner angle $\theta_j$
formed by the corresponding limiting tangents is well defined and satisfies
$\theta_j\in(0,2\pi)\setminus\{\pi\}$ for any $j=1,\ldots,m$;\\
(-) for each $p_j$ there exists a conformal bijection from an open neighborhood $U$ of $p_j$ which
maps $U\cap \pa\om$ onto a curve of class $C^2$.\\
In particular any regular domain is by definition simply connected.
\edf

We will use this definitions throughout the rest of this paper without further comments. Of course
polygons of any kind are regular according to our definition. The notations
$|\om|$ or $A(\om)$ will be used to denote the area of a simple domain $\om$,
while $L(\partial\om)$ will denote the length of the boundary of $\om$.

\brm\label{solregdef}
We will discuss at length solutions of a Liouville-type semilinear equation
with Dirichlet boundary conditions, see $P(\lm,\om)$ in section \ref{ss1.1} below.
In this respect, and if $\om$ is regular,
a solution $u$ will be by definition an $H^{1}_0(\om)$ weak solution \cite{GT} of the problem at hand,
$H^1_0(\om)$ being the closure of $C^{1}_c(\om)$ in the norm $\|u\|_2+\|\, |\nabla u|\,\|_2$.
In those cases where $\om$ is just assumed to be simple, a solution will be by definition a classical solution
$u\in C^{2}(\om)\cap C^{0}(\ov{\om})$.\\
It turns out that, by using the well known Brezis-Merle results \cite{bm} together with Lemma 2.1 in \cite{CCL},
any $H^{1}_0(\om)$ weak solution on a regular domain is also a classical
$C^{2}(\om)\cap C^{0}(\ov{\om})$ solution.
\erm

\bigskip

Let $\om\subset \R^2$ be open, bounded and simple. We define
$$
\mathcal{P}=\left\{\de\in L^{1}(\om)\,|\,\de\geq 0\;\mbox{a.e. in}\;\om,\;\ino\de =1\right\},
$$

and $G_\om(x,y)$ to be the unique  solution of

\beq\label{Green}
\graf{
-\Delta G_\om(x,y)= \delta_{x=y}&
\mbox{in}\hspace{.4cm}\;\; \om, \\
\hspace{0.7cm}G_\om(x,y)=0 &\mbox{on}\hspace{.1cm}\;\; \pa \om,}
\eeq
where $\delta_{x=y}$ is the Dirac distribution with singular point $y\in\om$,
$G_\om(x,y)=-\frac{1}{2\pi}\log(|x-y|)+H_\om(x,y)$ and $H_\om$ denotes the regular part.

For any $\de\in \mathcal{P}$ we also define the entropy and energy of $\de$ as
$$
\mathcal{S}(\de)=\ino\s(\de),\qquad
\mathcal{E}(\de)=\frac{1}{2}\ino \de G[\de],
$$
respectively, where
$$
\s(t)=\graf{-t\log{t}, &t>0,\\ 0, &t=0, \ }
$$
and
$$
G[\de](x)=\ino G_\om(x,y)\de(y)\,dy.
$$
For any $E\in\R$ we consider the MVP (Microcanonical Variational Principle)
$$
S(E)=\sup\left\{\mathcal{S}(\de),\quad \de\in \mathcal{P}_E \right\},\quad
\mathcal{P}_E=\{\de\in\mathcal{P}\,|\,\mathcal{E}(\de)=E\}.\qquad\qquad {\rm (MVP)}
$$

\bigskip

The following results has been obtained in \cite{clmp2} (see Propositions 2.1, 2.2, 2.3 in \cite{clmp2}):\\

MVP-(i) For any $E>0$, $S(E)<+\i$ and there exists $\de\in\mathcal{P}_E$ such that
$S(E)=\mathcal{S}(\de)$;\\
MVP-(ii) Let $\Upsilon=\frac{1}{|\om|}$ be the uniform density on $\om$ and
$E_\Upsilon=\mathcal{E}(\Upsilon)$. Then $\Upsilon$ is
a maximizer of $\mathcal{S}$ on $\mathcal{P}_{E_\Upsilon}$ and in particular if $|\om|=1$, then $S(E_\Upsilon)=0$; \\
MVP-(iii) If $|\om|=1$ then $S(E)$ is strictly increasing and negative for $E<E_\Upsilon$ and strictly decreasing and
negative for $E>E_\Upsilon$;\\
MVP-(iv) Let $\de^{\scp(E)}$ be a solution for the MVP at energy $E$. Then there exists $\bt=\bt_E\in\R$ such that
$$
\de^{\scp(E)}=\frac{e^{-\bt G[{\tsp \de}^{\scp(E)}]}}{\ino e^{-\bt G[{\tsp \de}^{\scp(E)}]}},
$$
or, equivalently, the function $\psi=G[\de^{\scp(E)}]$ satisfies the Mean Field Equation (MFE)
$$
\graf{
-\Delta \psi =\dfrac{e^{-\bt\psi}}{\ino e^{-\bt\psi}} & \mbox{in}\quad \om\\
\psi =0 & \mbox{on}\quad \pa\om
}\qquad\qquad (\mbox{MFE});
$$
MVP-(v) $S(E)$ is continuous.

\bigskip
\bigskip

We find it appropriate at this point to continue our discussion by introducing some concepts as
in \cite{clmp2} but with the aid of a slightly different mathematical arguments based on some
results in \cite{bm},  \cite{yy}, \cite{ls} and in particular in \cite{CCL} which were not at hand at that time.\\
Since solutions of the (MFE) with fixed $\bt>-8\pi$ are unique not only if $\om$ is simple and smooth \cite{suz}
but also if $\om$ is regular (see \cite{CCL}), and by using the Brezis-Merle \cite{bm}
theory of Liouville-type equations (as later improved in \cite{ls} and then in \cite{yy}) and the boundary
estimates in \cite{CCL}, we can divide the set of regular domains (see definition \ref{defsimp}) in two classes,
first introduced in \cite{clmp2}:

\bdf\label{kind}
Let $\om$ be regular. We say that $\om$ is of {\bf first kind} if
the unique (at fixed $\bt>-8\pi$ \cite{suz}, \cite{CCL}) solution $\psi_\bt$ of the {\rm (MFE)} satisfies
\beq\label{blow-up}
 \de_{(\bt)}:=\dfrac{e^{-\bt\psi}}{\ino e^{-\bt\psi}}\rightharpoonup\dt_{x=p},\quad\mbox{as}\quad\bt\searrow(-8\pi)^{+},
\eeq
weakly in the sense of measures, for some $p\in\om$.\\
We say that $\om$ is of {\bf second kind} otherwise.
\edf
We will skip the discussion of the case $\bt>0$ since its
mathematical-physical description is well understood \cite{clmp2}.\\
Let $\mathcal{E}(\de_{(\bt)})$ be the energy of the unique solution of the (MFE) with $\bt\in(-8\pi,0]$.
By using known arguments based on the results in \cite{bm}, \cite{ls} and  \cite{CCL}, \cite{yy} it can be shown that
either $\psi_\bt$ is uniformly bounded for $\bt\in(-8\pi,0]$ or it must satisfy \rife{blow-up} and in this case
in particular $\mathcal{E}(\de_{(\bt)})\to+\i$ as $\bt\searrow (-8\pi)^+$. Here is crucial Lemma 2.1 in \cite{CCL}
which ensures that solutions are
uniformly bounded in a neighborhood of $\pa\om$ whenever $\om$ is regular.

\brm
As a consequence of
an argument which we introduce in Lemma \ref{lem:231112} below, we could extend this alternative
(either $\psi_\bt$ is bounded or the energy $\mathcal{E}(\de_{(\bt)})\to+\i$ as $\bt\searrow (-8\pi)^+$)
to the case where $\om$ is just simple, the only difference
in this case being that one would have to allow (in principle) $p\in\pa\om$ in \rife{blow-up}.
However we do not know of any result claiming uniqueness of solutions
of the {\rm (MFE)} with $\bt\in(-8\pi,0)$ under such weak regularity assumptions on $\om$.
\erm

As in \cite{clmp2} we need the following:

\bdf\label{ecrit}
We set $E_c=\mathcal{E}(\de_{(\bt)})\left.\right|_{\bt=(-8\pi)^{+}}$ if $\om$ is of second kind and
$E_c=+\i$ if $\om$ is of first kind.
\edf

It has been shown in \cite{clmp2} that $E_\Upsilon<E_c$ and that to each $E_\Upsilon<E<E_c$ there corresponds a unique $\de^{\scp(E)}$ which attains the
supremum in the MVP and in particular a unique $\bt=\bt(E)\in (-8\pi,0)$ such that the corresponding unique
solution $\psi_\bt$ of the (MFE) satisfies $\de_{(\bt(E))}\equiv\de^{\scp(E)}$ and attains the supremum in the associated
CVP (Canonical Variational Principle)
$$
f(\bt)=f_\om(\bt)=\sup\{\mathcal{F}_\bt(\de),\;\;\de \in \mathcal{P}\,|\, -\mathcal{S}(\de)<+\i\},\qquad\qquad{\rm (CVP)}
$$
where, for $\de \in \mathcal{P}$,
$$
\mathcal{F}_\bt(\de)=-\frac{1}{\bt}\,\mathcal{S}(\de)+\mathcal{E}(\de),
$$
is the free energy of $\de$. In particular it has been proved in \cite{clmp2} that
$\mathcal{E}(\de_{(\bt)})$ is continuous and decreasing in $(-8\pi,0)$ and $S(E)$ is smooth and concave in $(E_\Upsilon,E_c)$.
Concerning these remarkable results
we refer to Theorem 3.1 and Proposition 3.3 in \cite{clmp2}.\\
In particular for domains of first kind the (mean field) thermodynamics of the system is rigorously defined for
any attainable value of the
energy and equivalently described by solutions of either the MVP or the CVP.  Actually, this problem is closely
related with another very subtle issue, that is, the fact
that solutions of the (MFE) always exist for $\bt\in (-8\pi,0]$  (a consequence of the Moser-Trudinger inequality
\cite{moser}) while in general do not exist for $\bt\leq -8\pi$, the value $\bt=-8\pi$ being the critical threshold
where the coercivity of the corresponding variational functional (that is \rife{var-func} below) breaks down.
A detailed discussion
of this point is behind our scopes and we
limit ourselves here with few details needed in the presentation of our results, see also section \ref{ss1.1} below.\\
Some sufficient conditions for the existence of solutions of the (MFE) at $\bt=-8\pi$ where provided in \cite{clmp1}
and hence used to show that for example
any long and thin enough rectangle is of second kind. The problem has been later solved in \cite{CCL} by
using a refined version of the subtle estimates in \cite{CLin1}, \cite{CLin2} and the newly derived uniqueness
of solutions of the (MFE) with $\bt\in(-8\pi,0]$ and, whenever they exist, for $\bt=-8\pi$ as well on regular domains.
In particular, it has been proved in Proposition 6.1 in \cite{CCL} that if $\om$ is regular, then
the following facts are equivalent:\\
SK-(i) $\om$ is of second kind;\\
SK-(ii) There is a solution of the (MFE) with $\bt=-8\pi$, say $\psi_{-8\pi}$;\\
SK-(iii) The unique branch of solutions of the (MFE) $\psi_\bt$ with $\bt\in(-8\pi,0]$ is uniformly bounded and
converges uniformly to $\psi_{-8\pi}$ as $\bt\searrow (-8\pi)^+$.\\

We conclude in particular that if the branch of (unique) maximizers satisfies \rife{blow-up},
then there is no solution of the (MFE)
with $\bt=-8\pi$ and in particular that a solution of the (MFE) with $\bt=-8\pi$ exists (and is unique) if and only
if {\it blow up for the {\rm (MFE)} at $\bt=8\pi$ occurs from the left}, that is, \rife{blow-up} occurs
but with $\bt\to(-8\pi)^{-}$. The fact that
(irrespective on the "side" which $\bt$ may choose to approach $8\pi$)
there is a branch of solutions which satisfy to a concentration property as in \rife{blow-up}, was already proved in
\cite{clmp2}, see NEQ-(ii) below.\\

The full theory as exposed in \cite{CCL} as well as the equivalence of statistical ensembles
has been recently extended to cover the case where $\om$ is multiply connected in \cite{BLin3}. As far as
one is concerned with the analytical problem of the existence for $\bt=-8\pi$ and uniqueness for $\bt\in[-8\pi,0)$, the
results in \cite{CCL} has been generalized in \cite{bl}, \cite{BLin2} to the case where
Dirac-type singular data are added in the (MFE).\\

\bigskip
The mean field thermodynamics for domains of second kind when $E\geq E_c$ is more involved.\\
Since it is not difficult to show that $\mathcal{F}_\bt$ is unbounded from above for $\bt<-8\pi$, then
there is no solution for the CVP with $\bt< -8\pi$ and therefore no equivalence (at all) among the MVP
and the CVP is at hand in this case. Nevertheless some insight about the range of energies $E\geq E_c$
was also obtained in \cite{clmp2}.
Let $\om$ be a domain of second kind. Then we have (see Propositions 6.1, 6.2 and Theorem 6.1 in \cite{clmp2}):\\
NEQ-(i) It holds
$$
-8\pi E+C_1 \leq S(E)\leq -8\pi E + C_2,\;\forall \,E\geq E_c,
$$
where $C_2=S(E_c)+8\pi E_c=8\pi f(-8\pi)$;\\
NEQ-(ii) Let $\de^{\scp (E)}$ be a solution of MVP at energy $E$. Then (up to subsequences)
$\de^{\scp (E)}\rightharpoonup\dt_{x=p}$, as $E\to+\i$, where $p$ is a maximum point of $H_\om(x,x)$;\\
NEQ-(iii) $S(E)$ is not concave for $E>E_c$.

\bigskip

Besides these facts, we do not know of any positive result about this problem for domains of second kind
when $E\geq E_c$.\\
It is one of our motivations to begin here a systematic study of the statistical mechanics
description of the case $E\geq E_c$. In this paper we work out the following program:\\
(-) Prove the existence of solutions of the (MFE) for suitable $\bt<-8\pi$ by assuming the domain
to be "thin" enough, see \S \ref{ss1.1} and \S \ref{ss1.4}.\\
(-) Prove that the first eigenvalue of the linearized problem for the (MFE) on those solutions
is strictly positive. This fact will imply that our solutions are local maximizers of ${\mathcal F}_\be$
as well as a multiplicity result yielding another set of unstable solutions, see \S \ref{ss1.2}.\\
(-) Prove that if the domain is "thin" enough, then there exists at most one solution of the (MFE) with
$\bt$ bounded from below and whose
energy is less than a certain threshold. This fact will imply that we
have found a connected and smooth branch of solutions where the energy is well defined as a
function of $\lm:=-\bt$, see Remark \ref{newrmsmooth} and \S \ref{ss1.3}.\\
(-) Prove that if the domain is "thin" enough and in a small enough range of energies, then the energy is
monotonic increasing as a function of $\lm=-\bt$. This fact will eventually imply that there exists one and only
one solution of the MFE at fixed energy (in that small range) which therefore is also the unique
maximizer of the entropy for the MVP. In particular we will prove that the entropy is concave in this range,
see \S \ref{ss1.4}.\\

\bigskip

This is the underlying idea which will guide us in the analysis of various problems of independent mathematical
interest as discussed in the rest of this introduction. We take the occasion here to provide all the motivations and/or
necessary comments about the statements
of the many results obtained (with the unique exception of Proposition \ref{pr2} below)
which is why the introduction is so lengthy.

\subsection{Existence of solutions for the supercritical (MFE) on thin domains}\label{ss1.1}$\left.\right.$\\
Amongst other things which will be discussed below, one of the main reasons which makes things
more difficult in the case $E\geq  E_c$ is the lack of a description of the solutions set for the
(MFE) with $\bt<-8\pi$. Since this will be a major point in our discussion, we introduce the quantities
$$
\lm:=-\bt,\qquad\mbox{and}\quad u=-\bt\psi=\lm\psi,
$$
and consider the following alternative but equivalent formulation of the (MFE)

$$
\graf{
-\Delta u =\lm \dfrac{e^u}{\ino e^u} & \mbox{in}\quad \om\\
u =0 & \mbox{on}\quad \pa\om
}\qquad P(\lm,\om)
$$

which we will denote by $P(\lm,\om)$. The following remark will be used throughout the rest of this paper.
\brm\label{rem3.1}
Clearly $P(\lm,\om)$ is rotational and translational invariant. Moreover
the integral in the denominator of the nonlinear datum in $P(\lm,\om)$ makes the problem \un{dilation invariant} too,
that is, $u$ is a solution of $P(\lm,\om)$ if and only if $v(y)=u(y_0+d_0R_0y)$ is a solution of
$P(\lm,\om^{(0)})$, where $y_0\in\R^2$, $d_0>0$, $R_0$ is an orthogonal $2\times2$ matrix and
$$
\om^{(0)}:=\{y\in\R^2\,|\,y_0+d_0R_0y\in\om\}.
$$
In particular, $u$ solves $P(\lm,\om_\a)$ with $\a=\frac{\ai}{\bi}$
where
\beq\label{1}
\om_\a=\{(x,y)\in \rdue\,|\,\a^2x^2+y^2\leq 1,\;\a\in (0,1]\},
\eeq
is the canonical two dimensional ellipse whose axis lengths are $\frac{1}{\a}$ and $1$,
if and only if $u_0(x^{'},y^{'})$ with $\{\bi x^{'}=\, x,\;\bi y^{'}= y\}$ solves $P(\lm,\mathbb{E}_{\ai,\bi})$, where
$$
\mathbb{E}_{\ai,\bi}=\{(x^{'},y^{'})\in \rdue\,|\,\ai^2 {x^{'}}^2+\bi^2{y^{'}}^2\leq 1,\;\ai\in (0,1],
\,\bi\in (0,1],\,\bi\geq \ai\},
$$
is the canonical two dimensional ellipse whose axis lengths are $\frac{1}{\ai}$ and $\frac{1}{\bi}$.
\erm

As mentioned above, we just miss a description of the solutions
set of $P(\lm,\om)$ with $\lm> 8\pi$ and $\om$ regular. General existence results for $P(\lm,\om)$ are at hand
for $\lm\in \R\setminus{8\pi\enne}$ only if $\om$ is a multiply connected domain,
see \cite{DJLW}, \cite{st} and the deep results in \cite{CLin2} (see also \cite{BDeM}).\\
This is far from being a technical problem. Indeed, a well known result based on the Pohozaev identity (see for example
\cite{clmp1}) shows that if $\om$ is strictly starshaped, then there exists
$\lm_*=\lm_*(\om)\geq 8\pi$ (see also Remark \ref{extremal} below) such that $P(\lm,\om)$ has no solutions for
$\lm\geq \lm_*(\om)$.
This result is sharp since indeed
$\lm_*(B_R(0))=8\pi$, where $B_R(0)=\{x\in\R^2\,:\,|x|<R\}$ for some $R>0$.\\
Therefore, in particular, the Leray-Shauder degree
of the resolvent operator for $P(\lm,\om)$ with $\om$ regular vanishes identically for any
$\lm>8\pi$, see \cite{CLin2}.\\
If this were not enough we also observe that, at least in case $\om$ is convex,
the well known results in \cite{BaPa}, \cite{CLin1}, \cite{EGP}, \cite{KMdP} concerning concentrating solutions for
$P(\lm,\om)$ as $\lm\to 8\pi k$, for some fixed $k\in \enne$, are of no help here,
since it has been shown in \cite{GrT} that in fact neither those blow-up solutions sequences
exist if $k\geq 2$.\\
Finally let us remark that we are concerned here just with solutions of $P(\lm,\om)$. If
we allow some weight to multiply the exponential nonlinearity, then other solutions exist for $\lm>8\pi$ on simply
connected domains, see
for example \cite{B1-1}, \cite{B2}, \cite{BM2} and more recently the general results derived in \cite{BMal}.\\
As a matter of fact, the only general result we are left with is the immediate corollary of the uniqueness
results in \cite{CCL}, which shows that:\\

SK-(iv) if $\om$ is of second kind, then the branch of unique solutions
$u_{\scp \lm}$, $\lm\in[0,8\pi]$ of $P(\lm,\om)$ can be extended (via the implicit function theorem)
in a small right neighborhood of $8\pi$.\\

Our first result is concerned with a sufficient condition for the existence of solutions of $P(\lm,\om)$ with
$\lm>8\pi$ on "thin" domains.

\bte\label{t1-intro}$\left.\right.$\\
{\rm (a)} Let $\om$ be a simple domain.
For any $c\in(0,1]$ there exist $\overline{\a}_*>\underline{\a}_*(c)>0$ such that
if $\{\a^2x^2+y^2\leq \beta_-^2\}\subset\Omega\subset\{\a^2x^2+y^2\leq\beta_+^2\}$
with $c=\tfrac{\beta_-^2}{\beta_+^2}$ then, for
any $\a\in(0,\underline{\a}_*(c)]$ and for any $\lm\leq\lm_{\a,c}$, there exists a solution $\ul$ of
$P(\lm,\om)$, where $\underline{\lm}_{\a,c}< \lm_{\a,c}<\overline{\lm}_\a$ and $\underline{\lm}_{\a,c}$, $\overline{\lm}_\a$ are strictly decreasing (as functions of $\rho$) in $(0,\underline{\a}_*(c)]$, $(0,\overline{\a}_*]$ respectively with
$\underline{\lm}_{\underline{\a}_*(c),c}=8\pi=\overline{\lm}_{\overline{\a}_*}$ and $\underline{\lm}_{\a,c}\simeq
\frac{4\pi c}{(8-c)\a}$, $\ov{\lm}_{\a}\simeq
\frac{11\pi}{16\a}$ as $\a\to 0^+$.\\

{\rm (b)}
There exists $\bar N> 4\pi$ such that if $\Omega$ is an open, bounded and convex set (therefore simple)
whose isoperimetric ratio, $N\equiv N(\Omega)=\frac{L^2(\pa\Omega)}{A(\Omega)}$, satisfies $N\geq\bar N$, then
for any $\lm\leq\lm_\textnormal{\tiny{$N$}}$ there exists a
solution $\ul$ of $P(\lm,\Omega)$, where
$\underline{\Lambda}_\textnormal{\tiny{$N$}}<\lambda_\textnormal{\tiny{$N$}}<
\overline{\Lambda}_\textnormal{\tiny{$N$}}$ with $\underline{\Lambda}_\textnormal{\tiny{$\bar N$}}=
8\pi$, $\underline{\Lambda}_\textnormal{\tiny{$N$}}$ and $\overline{\Lambda}_\textnormal{\tiny{$N$}}$
strictly increasing in $N$ and $\underline{\Lambda}_\textnormal{\tiny{$N$}}\simeq\frac{\pi^2N}{496}+O(1)$,
$\overline{\Lambda}_\textnormal{\tiny{$N$}}\simeq\frac{33\sqrt3N}{16\pi}+O(1)$ as $N\to+\infty$.
\ete

\brm\label{susp}
The suspect that this result should hold was initially due to the above mentioned result in \cite{clmp1}
(which states that if
$\om$ is a long and thin enough rectangle then a solution of $P(8\pi,\om)$ exists) and to a result in \cite{CCL}
(which states that there exists a critical value $d_1<1$ such that if $\om$ is a rectangle whose sides lengths are
$a_1\leq b_1$, then a solution of $P(8\pi,\om)$ exists if and only if $\frac{a_1}{b_1}\leq d_1$). In particular this
observation already shows that ${\bar N}>4\pi$.
\erm

\brm\label{extremal}
Clearly $c=1$ if and only if $\Omega$ is an ellipse,
while if $\Omega$ is a rectangle it is easy to see that $c=\frac12$ is optimal.
We also
have the quantitative estimate $0.0702<\overline{\a}_*(1)$ which could be used in principle
to obtain an estimate for either $d_1$ (see Remark \ref{susp}) or $\bar N$. We will not insist about this point since it seems that we are too
far from optimality.
In the case of the ellipse $\om_\a$, the existence lower/upper threshold values
$\un{\lm}_\a\simeq\frac{4\pi}{7\a}/\ov{\lm}_{\a}\simeq
\frac{11\pi}{16\a}$ should
be compared with the Pohozaev's upper bound for the existence of solutions for $P(\lm,\om_\a)$, that is
$$
\lm <\lm_*(\om_\a):=4\int\limits_{\pa \om_\a}\frac{ds}{(\underline{x},\underline{\nu}\,)}=\frac{4\pi}{\a}(1+\a^2).
$$
\erm

\brm\label{rem:branch}
For regular domains, the branches of solutions obtained above will be seen to be connected and smooth, see Remark
\ref{newrmsmooth} below. We will denote them by $\mathcal{G}_{\a,c}=\{(\lm,\ul)\,:\,\lm\in[0,\lm_{\a,c}]\}$ (as obtained in Theorem
\ref{t1-intro}(a)) and
$\mathcal{G}_{N}=\{(\lm,\ul)\,:\,\lm\in[0,\lm_\textnormal{\tiny{$N$}}]\}$ (as obtained in Theorem
\ref{t1-intro}(b)) respectively.
\erm

The proof of Theorem \ref{t1-intro} is, surprisingly enough, based on the sub-supersolutions method. In particular
we use the result in \cite{clsw} which allows for such a weak assumptions about the regularity of $\om$.
The underlying idea in case $\om=\om_\a$ is:\\
(-) if the ellipse $\om=\om_\a$ is "thin" enough (i.e. if $\a$ is small enough)
then the branch of \un{minimal} solutions for the classical Liouville problem
$$
\graf{
-\Delta u=\mu \,{\dsp e^u} & \mbox{in}\quad \om\\
u=0 & \mbox{on}\quad \pa\om
}\qquad Q(\mu,\om)
$$
cannot be pointwise too far from the $C^{2}_0(\om_\a)$ function
$$
v_{\a,\ga}=2\log{\left(\frac{1+\ga^2}{1+\ga^2(\a^2x^2+y^2)}\right)},\quad (x,y)\in \om_\a,
$$
for a suitable value of $\ga$ depending on $\mu$ and $\a$. Of course, the guess about $v_{\a,\ga}$ is
inspired by the Liouville formula \cite{Lio}. Therefore, for fixed $\mu$ and $\a$, we seek values
$\ga_{\mp}$ such that
$v_{\a,\ga_{\mp}}$ are sub-supersolutions respectively of $Q(\mu,\om_\a)$.\\
(-) if the choice of $\ga_{\pm}(\mu)$ is made with enough care, then, along the branch of solutions
(say $u_{\scp \mu}$) for $Q(\mu,\om)$ found via the sub-supersolutions method,
the value of $\lm$ \un{defined} as follows
$$
\lm:=\mu\int\limits_{\om_\a} e^{u_{\scp \mu}},
$$
can be quite large whenever $\a$ is small enough.

\bigskip

Part (b) of Theorem \ref{t1-intro} will be a consequence of Part (a) and Theorems \ref{tjohn} and \ref{tlassek}
below.

\bte\label{tjohn}{{\rm \{}\cite{John}{\rm \}}}
Let $K\subset\R^2$ be a convex body (that is a compact convex set with nonempty interior).
Then there is an ellipsoid $E$ (called the John ellipsoid which is the ellipsoid of maximal volume contained in $K$)
such that, if $c_0$ is the center of $E$, then the inclusions
$$
E\subset K\subset \{c_0+2(x-c_0)\,:\,x\in E\}
$$
hold.
\ete

\bte\label{tlassek}{{\rm \{}\cite{Lassek-priv}{\rm \}}}
Every convex body $K\subset\R^2$ contains an ellipse of area $\tfrac{\pi}{3\sqrt3}\,A(K)$.
\ete

A short proof of the previous theorem is based on a result in \cite{Besicovitch}, where the existence of an
affine-regular hexagon $H$ of area at least $\tfrac23\,A(K)$ and inscribed in $K$
is established. Indeed, considering the concentric inscribed ellipse in $H$ one gets the thesis.

\brm\label{rem-Lassek}
In particular Theorem \ref{tlassek} has been used to obtain the asymptotic behaviors of
$\underline{\Lambda}_\textnormal{\tiny{$N$}}$ and $\overline{\Lambda}_\textnormal{\tiny{$N$}}$.
A more rough estimate of those asymptotics could have been obtained by using other (much worst)
known estimates of the area of the enclosed ellipse. In particular, while Theorem \ref{tjohn}
is well known \cite{John}, it seems that Theorem \ref{tlassek} is not
and we are indebted with Prof. M. Lassak who kindly reported to us a proof of it \cite{Lassek-priv}
based on the cited reference \cite{Besicovitch}.
\erm

\bigskip

Clearly, as an immediate corollary of Theorem \ref{t1-intro} and the equivalence of
SK-(i) and SK-(ii) we conclude that
if $\om$ is regular and satisfies the assumptions of Theorem \ref{t1-intro}(a)(Theorem \ref{t1-intro}(b))
with $\a\in(0,\un{\a}_*(c)]$ ($N(\om)>\bar N$) then it is of second kind.

\subsection{Non degeneracy and multiplicity of solutions of the supercritical (MFE) on thin domains}\label{ss1.2}$\left.\right.$\\
Let us define the density corresponding to a solution $\ull$ of $P(\lm,\om)$ as
\beq\label{dedef}
\de_{\scp \lm}\equiv\de(\ull):=\dfrac{\e{\ull}}{\int\limits_{\om} \e{\ull}}.
\eeq
A crucial tool used in the proof of the equivalence of statistical ensembles \cite{clmp2}
is the uniqueness of solutions \cite{suz}, \cite{CCL} (see also \cite{BLin3})
of $P(\lm,\om)$ for $\lm\in[0,8\pi]$.
The situation is far more involved in case $\lm>8\pi$ since on domains of second kind, solutions
are not anymore unique.\\
This fact is already clear from NEQ-(ii) and SK-(iv) above, that is, if $\om$ is of second kind
we have a blow-up branch which satisfies
\beq\label{blow-up2}
\de(\ull)\rightharpoonup\dt_{x=p},\quad\mbox{as}\quad\lm\searrow(8\pi)^{+},
\eeq
weakly in the sense of measures, for some critical point $p\in\om$ of $H_\om(x,x)$, and the smooth solutions
of $P(\lm,\om)$ in a small right neighborhood of $8\pi$. Hence, we have at least two solutions in a right
neighborhood of $8\pi$, a well known fact that could have been also deduced
by using the alternative in Theorem 7.1 in \cite{clmp2} together with the uniqueness result in \cite{CCL}.\\
We wish to make a further step in this direction.
To this purpose we first study the linearized problem for $P(\lambda,\Omega)$ at $\ul$,
where $\ul$ is the solution obtained in Theorem \ref{t1-intro},
showing the positivity of its first eigenvalue (see Proposition \ref{pr2} and Remark \ref{branch} for details).
It is worth to point out that the above fact, which yields a multiplicity result too,
is also crucial in the analysis of the solutions branches $\mathcal{G}_{\a,c},\mathcal{G}_N$, see Remarks
\ref{rem:branch} and \ref{newrmsmooth}. In particular we have:

\bpr\label{pr3} For fixed $c\in(0,1]$, let $\Omega$ be a regular domain that
satisfies $\{\a^2x^2+y^2\leq \beta_-^2\}\subset\Omega\subset\{\a^2x^2+y^2\leq\beta_+^2\}$,
with $\frac{\bt^2_-}{\bt^2_+}=c$ and $\a\in(0,\underline{\a}_*(c)]$,
with $\underline{\a}_*(c)$ as found in Theorem \ref{t1-intro}(a).
Let $\om$ be a convex domain with $N(\om)>\bar N$ as found in Theorem \ref{t1-intro}(b).\\
The portions of $\mathcal{G}_{\a,c},\mathcal{G}_{\textnormal{\tiny{$N$}}}$ with $\lm\in[0,8\pi]$ coincide
with the branch of unique absolute minimizers of
\beq\label{var-func}
F_{\lm}(u)=\frac12\int\limits_{\om} |\nabla u|^2 \,dx-\lm
\log\left(\;\int\limits_{\om} e^u\,dx\;\right),\quad u\in H^{1}_0(\om),
\eeq
and for each $\lm\in(8\pi,\lm_{\a,c}]$ or $\lm\in(8\pi,\lm_{\textnormal{\tiny{$N$}}}]$ the corresponding
solutions $\ul$ such that $(\lm,\ul)\in \mathcal{G}_{\a,c}$ and $(\lm,\ul)\in \mathcal{G}_\textnormal{\tiny{$N$}}$
are strict local minimizers of $F_{\lm}$.
\epr

\brm\label{newrmsmooth}
By using the bounds provided by the sub-supersolutions method (see \rife{susu} in the proof of Theorem \ref{t1-intro}), 
Proposition \ref{pr2} and Theorem \ref{unique:1-intro} below, then standard bifurcation theory \cite{CrRab}
shows that
for any fixed $\ov{\lm}>8\pi$, possibly taking a smaller $\un{\a}_*(c)$ and a larger $N$,
the portions of $\mathcal{G}_{\a,c}$ and $\mathcal{G}_\textnormal{\tiny{$N$}}$ with
$\lm\leq \ov{\lm}$ are smooth and connected branches with no bifurcation points.
\erm

The proof of Proposition \ref{pr3} is a straightforward consequence of the fact that the first eigenvalue of the
linearized problem for $P(\lm,\om)$ is strictly positive along $\mathcal{G}_{\a,c}$ and
$\mathcal{G}_{\textnormal{\tiny{$N$}}}$, see Proposition \ref{pr2} in section \ref{sec2}.\\
We shall see that, by virtue of Proposition \ref{pr3},
it is possible to show that for $\lm\in (8\pi,\lm_{\a,c})\setminus 8\pi\enne$
the functional $F_{\lm}$ exhibits a mountain-pass type structure which in turn yields
the existence of min-max type solutions to $P(\lm,\Omega)$. More precisely we obtain the following result.

\bte\label{mp-intro}$\left.\right.$\\
{\rm (a)} Let $\om$, $\a\in(0,\un{\a}_*(c)]$ and $\lm_{\a,c}$ be as in Theorem \ref{t1-intro}(a)
and let $\ul$ be a solution of $P(\lm,\om)$ for $\lm\leq\lm_{\a,c}$.
Then, for any $\lm\in(8\pi,\lm_{\a,c})\setminus 8\pi\enne$ there exists a second solution
$\vl$ of $P(\lm,\Omega)$ such that $F_{\lm}(\vl)> F_{\lm}(\ul)$. \\
{\rm (b)} Let $\om$, $\bar N>4\pi$, $N(\om)$ and $\lm_\textnormal{\tiny{$N$}}$ be as
in Theorem \ref{t1-intro}(b)
and let $\ul$ be a solution of $P(\lm,\om)$ for $\lm\leq\lm_\textnormal{\tiny{$N$}}$.
Then, for any $\lm\in(8\pi,\lm_\textnormal{\tiny{$N$}})\setminus 8\pi\enne$ there exists a second solution
$\vl$ of $P(\lm,\Omega)$ such that
$F_{\lm}(\vl)> F_{\lm}(\ul)$.
\ete

\brm By using well known compactness results \cite{yy} as well as those recently derived in \cite{GrT}, we
conclude that any sequence of solutions $\vl$ with $8\pi k<\lm< 8\pi(k+1),\,k\geq 1$ obtained in part
{\rm (b)} converges as $\lm\to 8\pi (k+1)$ to a solution $v_{8\pi(k+1)}$ of $P(8\pi(k+1),\om)$.
We also have at least two different arguments showing that for any fixed $\ov{\lm}>0$,
possibly taking  a larger $N,$ those $v_{8\pi k}$
which also satisfy $8\pi k\leq\ov{\lm}$ are distinct from
those obtained in Theorem \ref{t1-intro}{\rm (b)} for $\lm=8\pi k$.
The first one is a standard bifurcation-type argument
based on Remark \ref{newrmsmooth} and Proposition \ref{pr2} below. The second one is based
on the uniqueness result stated in Theorem \ref{unique:1-intro} below. 
\erm

\brm\label{metastable}
It is easy to check that if $u$ is a solution of $P(\lm,\om)$ and $\de(u)$ is defined as in \rife{dedef}, then
$\de(u)$ is a critical point of $\mathcal{F}_{-\lm}$ and in particular
$\mathcal{F}_{-\lm}(\de)=-\frac{1}{\lm^2}F_\lm(u)$. Hence, if $\ul$ and $\vl$ are as in Theorem \ref{mp-intro}
with $\de(\ul)$ and $\de(\vl)$ as in \rife{dedef}, then it is readily seen that
$\mathcal{F}_{-\lm}(\de(\ul))<\mathcal{F}_{-\lm}(\de(\vl))$. In particular $\de(\ul)$ is a kind
of metastable state
(in the sense that it is a strict local maximizer of $\mathcal{F}_{-\lm}$) while $\de(\vl)$ is expected to be
unstable (since it is a min-max type critical point of $\mathcal{F}_{-\lm}$).\\
In any case, whenever $\om$ is regular (and since solutions of $P(8\pi,\om)$ are unique in this case \cite{CCL}),
then any sequence of solutions found in Theorem \ref{mp-intro} for $P(\lm,\om)$ with $\lm\searrow 8\pi^+$
must satisfy \rife{blow-up2}.

\erm

\subsection{Uniqueness of solutions for the supercritical (MFE) with bounded energy on thin domains}\label{ss1.3}$\left.\right.$\\

As a matter of fact we are still unable to define the energy as a monodrome function
of $\lm$. We explain the next step toward this goal in the case of the ellipse $\om_\a$.\\

Although solutions of $P(\lm,\om_\a)$ are not unique as a function of $\lm$, what we can prove is that
for fixed $\ov{\lm}\geq 8\pi$ and $\ov{E}\geq 1$, then for $\a$ small enough there could be at most
one solution $u_{\scp \a, \lm}$ such that $\lm\leq\ov{\lm}$ and
\beq\label{endef2-intro}
\mathcal{E}(\de(u_{\scp \a, \lm}))\leq \ov{E}.
\eeq

This is a major achievement since, by using also Proposition \ref{pr2} below, it implies that
(as far as $\a$ is small enough),
the energy (see Proposition \ref{pr3:I}) is well defined as a function of $\lm$, whenever $\lm\leq \ov{\lm}$
and the supremum of the range of the energy itself is not greater than $\ov{E}$.\\
Let us think at the results obtained in \S \ref{ss1.1} and \S \ref{ss1.2} in terms of the
$(\lm,\|u_{\scp \a, \lm}\|_\infty)$ bifurcation diagram. To fix the ideas,
we propose the following naive description. As $\a$ gets smaller and smaller, we have:\\
(-)  the portion with $\lm\leq\ov{\lm}$ and $\mathcal{E}(u_{\scp \a,\lm})\leq \ov{E}$
of the (smooth, see Remark \ref{newrmsmooth}) branches of solutions $\mathcal{G}_{\a,c},\mathcal{G}_{N}$ obtained in Theorem \ref{t1-intro}
gets lower and flatter, that is, $\|u_{\scp \a,\lm}\|_\infty\searrow 0^+$.
See also Remark \ref{rem6.1} below.\\
(-) In the same time the portion with $\lm\leq\ov{\lm}$ of the branches obtained in Theorem \ref{mp-intro}
(as well as any other possible solution) gets higher and higher the corresponding energies getting greater and finally greater than $\ov{E}$.\\
(-) Any bifurcation/bending point one should possibly meet along $\mathcal{G}_{\a,c},\mathcal{G}_{N}$
moves in the region $\lm>\ov{\lm}$.

\bigskip

It is understood that the value $1$ in the condition $\ov{E}\geq 1$ could have been substituted by any other
fixed positive number. More exactly we have the following:
\bte\label{unique:1-intro}$\left.\right.$\\
Fix $\ov{\lm}\geq 8\pi$ and $\ov{E}\geq 1$. Then:\\
{\rm (a)} Let $\om$ be a simple domain and suppose that there exists $c\in(0,1]$
such that $\{\a^2x^2+y^2\leq \beta_-^2\}\subseteq\Omega\subseteq\{\a^2x^2+y^2\leq\beta_+^2\}$
with $c=\tfrac{\beta_-^2}{\beta_+^2}$.\\ Then there exists
$\widetilde{\a}_{1}=\widetilde{\a}_{1}(c,\ov{E},\ov{\lm})>0$ such that for any $\a\in(0,\widetilde{\a}_{1}]$,
there exists at most one solution $u_{\scp \lm}$ of $P(\lm,\om)$ with $\lm\leq\ov{\lm}$ which satisfies
\rife{endef2-intro}.\\

{\rm (b)} Let $\om$ be any open, bounded and convex (therefore simple) domain.
There exists $\widetilde{N}=\widetilde{N}(\ov{\lm},\ov{E})\geq  4\pi$ such that for any such $\om$ satisfying
$$
N(\om):=\frac{L^2(\pa \om)}{A(\om)}\geq  \widetilde{N},
$$
there exists at most one solution $u_{\scp \lm}$ of $P(\lm,\om)$ with $\lm\leq\ov{\lm}$ which satisfies
\rife{endef2-intro}.
\ete

The proof of Theorem \ref{unique:1-intro} is based on two main tools.\\
The first one is
an a priori estimate for solutions of $P(\lm,\om)$ (which satisfy $\lm\leq\ov{\lm}$ and \rife{endef2-intro})
with a uniform constant $\ov{C}$ which do not depend neither on $u$ nor \un{on the domain $\om$}, but only on
$\ov{\lm}$ and $\ov{E}$. Roughly speaking,
and in case $\om=\om_\a$, this kind of uniformity with respect to the domain is needed since we consider the limit in which
$\a$ gets very small, that is, we seek uniqueness for \un{all} domains which are "thin" in the sense specified
in the statement of Theorem \ref{unique:1-intro}. We refer to
Lemma \ref{lem:231112} and the discussion about it in section \ref{sec:unique} for further details.\\
The second tool is a careful use of the dilation invariance (see Remark \ref{rem3.1}) to be used together with
an estimate about the first eigenvalue of the Laplace-Dirichlet problem on a "thin" domain,
see \rife{110213.21} below for more details.

\subsection{Uniqueness of solutions for the supercritical (MFE) on $\om_\a$ with fixed energy and concavity of the
Entropy }\label{ss1.4}$\left.\right.$\\
In this subsection we fix $\om=\om_\a$.\\
As observed above, by using Theorem \ref{unique:1-intro} and Proposition \ref{pr2} below
we can prove that (as far as $\a$ is small enough) the energy (see Proposition \ref{pr3:I}) is well
defined as a function of $\lm$ (along the branch $\mathcal{G}_{\a,1}$ found in
Theorem \ref{t1-intro}(a), see Remark \ref{newrmsmooth}) whenever $\lm\leq \ov{\lm}$ and the supremum
of the range of the energy itself
is not greater than $\ov{E}$. It is tempting at this point to say that the entropy maximizers of the
MVP are those solutions of the (MFE) obtained in Theorem \ref{t1-intro}(a). However we still
don't know whether or not this is true, since obviously there could be many solutions on $\mathcal{G}_{\a,1}$ (i.e.
with different values of $\lm$) corresponding to a fixed energy $E\leq \ov{E}$
(see for example fig. 5 in \cite{clmp2}).
In such a situation it would be difficult to detect which is, (or worst, which are) the one which
really maximizes the entropy.
A possible solution to this problem could be obtained if we would be able to understand the monotonicity
of the energy as a function of $\lm$ on $\mathcal{G}_{\a,1}$. The first step toward this goal is to show
that the solutions of $P(\lm,\om_\a)$ obtained in Theorem \ref{t1-intro}(a) can be
expanded in powers of $\a$ with the leading order taking up an explicit and simple form (see also
\rife{191112.1}, \rife{021212.2} below), that is
\beq\label{phi0-intro}
\phi_0(x,y;\lm,\a)=\mu_0(\lm,\a)\psi_0(x,y;\a),\quad (x,y)\in \om_\a,
\eeq
where $\mu_0$ satisfies \rife{mu0-intro}-\rife{mu0-intro1} below and
\beq\label{080213.1-intro}
\psi_0(x,y;\a)=\frac{1}{2(1+\a^2)}\left(1-(\a^2x^2+y^2)\right),\quad (x,y)\in \om_\a.
\eeq
Of course, we could have used the fact that we already knew about the existence of the branch
$\mathcal{G}_{\a,1}$ and managed to expand those solutions as a function of $\a$.
Instead we decided to make the argument self-contained by pursuing
another proof of independent interest of the existence of solutions of $P(\lm,\om_\a)$. It shows that there exists
$\a_0$ small enough (depending on $\ov{\lm}$) such that for any $\a<\a_0$ and for each $\lm\in[0,\ov{\lm})$
a solution $\ull$ for $P(\lm,\om_\a)$ exists whose leading order with respect to $\a$
takes up the form \rife{phi0-intro}. There is no problem in checking that these solutions coincide with
those on the branch $\mathcal{G}_{\a,1}$ obtained in Theorem \ref{t1-intro}(a). Indeed this is at this point
an easy consequence of Theorem \ref{unique:1-intro}.\\
We still face the problem of how to handle the term $\int\limits_{\om_\a}e^{\ull}$
in the denominator of the nonlinear term in $P(\lm,\om_\a)$. This time we will solve this issue by
seeking solutions $v_\a$ of $Q(\mu_0\rho,\om_\a)$ which satisfy the following identity in a suitable
set of values of $\lm$,
\beq\label{lm0-intro}
\lm=\mu_0\a\int\limits_{\om_\a}e^{\ull}.
\eeq

This is the content of Theorem \ref{thm:261112-intro} below. More exactly, by setting
$$
D^{(k)}_\lm=\frac{\pa^{k}}{\pa\lm^k},\;k=0,1,2,
$$
we have the following:
\bte\label{thm:261112-intro}
Let $\ov{\lm}\geq 8\pi$ be fixed. There exists $\a_0>0$ depending on $\ov{\lm}$ such that for
any $\a<\a_0$ and for each $\lm\in[0,\ov{\lm}\,)$ there exists a solution $\ull$ for $P(\lm,\om_\a)$
which satisfies
\beq\label{sol:II-intro}
\ull(x,y;\lm)=\a\phi_0(x,y;\lm)+\a^2\phi_{1}(x,y;\lm)+\a^3\phi_2(x,y;\lm),
\quad (x,y)\in \om_\a,
\eeq
where $\{\phi_0,\phi_{1},\phi_2\}\subset C^{2}_0(\om)$. Moreover $\phi_0$ takes the form \rife{phi0-intro} with
$\mu_0$ a smooth function which satisfies
\beq\label{mu0-intro}
\mu_0(\lm,\a)=\frac{\lm}{\pi}-\frac{\lm^2}{4\pi^2}\a+\mbox{\rm O}(\a^2),
\eeq
and
\beq\label{mu0-intro1}
D^{(1)}_\lm \mu_0(\lm,\a)=\frac{1}{\pi}-\frac{\lm}{2\pi^2}\a+\mbox{\rm O}(\a^2),\quad
D^{(2)}_\lm \mu_0(\lm,\a)=-\frac{1}{2\pi^2}\a+\mbox{\rm O}(\a^2).
\eeq
In particular the following uniform estimates hold
\beq\label{regular-intro}
\|D^{(k)}_\lm\phi_0\|_{\scp C^{2}_0(\om)}+\|D^{(k)}_\lm\phi_{1}\|_{\scp C^{2}_0(\om)}+
\|D^{(k)}_\lm\phi_2\|_{\scp C^{2}_0(\om)}\leq \ov{M}_k,\;k=0,1,2,
\eeq
for suitable constants $\ov{M}_k$, $k=0,1,2$ depending only on $\ov{\lm}$.
Finally these solutions' set is a smooth branch which coincides with a portion of $\mathcal{G}_{\a,1}$.
\ete

\brm\label{rem:pr2}
In the proof of Theorem \ref{thm:261112-intro} and therefore
in \un{all} the expansions in powers of $\a$ what we really use is
the fact that solutions $v_\a$ of $Q(\mu_0\a,\om_\a)$ can be expanded in powers of $\a$ and in particular that
$\lm_0(\mu_0,\a):=\mu_0\a\int\limits_{\om_\a}e^{v_\a}$ is smooth, see Lemma \ref{intermedio} below.
Here we need some estimates about the first eigenvalue of the linearization of $Q(\mu,\om)$ as obtained in
Proposition \ref{pr2} below.
\erm

By using Theorem \ref{thm:261112-intro} we can prove the following result.
Let $\widetilde{\a}_{1}$ be fixed as in Theorem \ref{unique:1-intro}(a). Then we have:

\bte\label{unique:2-intro}
Let $\ov{\lm}\geq 8\pi$ and let $\widehat{E}_\a$ be defined by
$$
\widehat{E}_\a:=\frac{\a}{8\pi}+\frac{\a^2}{50\pi^2}\ov{\lm}.
$$
For each $\a<\widetilde{\a}_{1}$ and $E\in\left[\frac{\a}{8\pi},\widehat{E}_\a\right]$ there exists one and only one
solution $u_{\scp \lm}$
for $P(\lm,\om_\a)$ whose energy is $\mathcal{E}(\de(u_{\scp \lm}))=E$. Let $\widehat{\lm}_\a$ be defined by
$\mathcal{E}(\de(u_{\scp \widehat{\lm}_\a}))=\widehat{E}_\a$. Then in particular the identities
$$
\widehat{E}(\lm)=\mathcal{E}(\de(u_{\scp \lm})),\quad \mathcal{E}(\de(u_{\scp \widehat{\lm}(E)}))=E,
$$
define: \\
$\widehat{E}(\lm):[0,\widehat{\lm}_\a]\to\left[\frac{\a}{8\pi},\widehat{E}_\a\right]$ as a smooth and strictly
increasing function of $\lm$ and\\
$\widehat{\lm}(E):\left[\frac{\a}{8\pi},\widehat{E}_\a\right]\to[0,\widehat{\lm}_\a]$
as a smooth and strictly increasing function of $E$.\\

Moreover we have

\beq\label{231112.3-intro}
\widehat{E}(\lm)=\frac{\a}{8\pi}+\frac{\a^2}{48\pi^2}\lm+\mbox{\rm O}(\a^3),\quad
\widehat{\lm}(E)=\frac{48\pi^2}{\a^2}\left(E-\frac{\a}{8\pi}\right)+\mbox{\rm O}(\a).
\eeq
\beq\label{231112.3.a-intro}
\frac{d}{d\lm}\widehat{E}(\lm)=\frac{\a^2}{48\pi^2}+\mbox{\rm O}(\a^3),\quad
\frac{d}{d E}\widehat{\lm}(E)=\frac{48\pi^2}{\a^2}+\mbox{\rm O}(\a),
\eeq
\beq\label{231112.3.b-intro}
\frac{d^2}{d\lm^2}\widehat{E}(\lm)=\mbox{\rm O}(\a^3),\quad
\frac{d^2}{d E^2}\widehat{\lm}(E)=\mbox{\rm O}(\a).
\eeq

\ete

\brm\label{rem6.1}
The notation $\mbox{\rm O}(\a^m)$, $m\in\enne$ is used here and in the rest of this paper to denote various quantities
uniformly bounded by $C_m\a^m$ with $C_m>0$ a suitable constant depending only on $\ov{\lm}$.\\
This result is
consistent with the underlying idea that, as $\a$ gets smaller and
smaller, then the energies of the entropy maximizers (which are solutions of $P(\lm,\om_\a)$)
with values of $\lm$ uniformly bounded from
above have to approach the energy of the uniform density distribution
$\Upsilon=\frac{1}{|\om_\a|}$, that is
$$
E_{\Upsilon,\a}:=\mathcal{E}\left(\frac{1}{|\om_\a|}\right)=\frac12\int\limits_{\om_\a}
\frac{1}{|\om_\a|}G_\a\left[\frac{1}{|\om_\a|}\right]=
\frac{\a}{2\pi}\int\limits_{\om_\a}\frac{1}{|\om_\a|2(1+\a^2)}\left(1-(\a^2x^2+y^2)\right)=\frac{\a}{8\pi}.
$$
Here we used the easily derived explicit expression of the function
$G_\a\left[\frac{1}{|\om_\a|}\right]$ see also \rife{phi0-intro}, \rife{080213.1-intro} and \rife{191112.1},
\rife{021212.2} below.
\erm

\brm
In particular \rife{231112.3-intro} yields $\widehat{\lm}_\a=\frac{48}{50}\ov{\lm}+\mbox{\rm O}(\a)$ and since
$\ov{\lm}\geq 8\pi$ can be chosen at wish and (see Definition \ref{ecrit}) $E_c=\mathcal{E}(\de(u_{8\pi}))$, then of course $E_{\Upsilon,\a}<E_c<\widehat{E}_\a$ and
we succeed in the description of the energy as a function of (minus) the inverse temperature $\lm=-\be$ in
a very small range of energies above $E_c$.
\erm

Let us observe that \rife{231112.3.a-intro} is in perfect agreement with the discussion in \S \ref{ss1.3},
that is, the portion with $\lm\leq\ov{\lm}$
of the branch of solutions obtained in Theorem \ref{t1-intro} gets lower and flatter as $\a$ gets smaller and smaller.
Actually we could not find another way to prove Theorem \ref{unique:2-intro} than explicit evaluations.
This is why our concern in Theorem \ref{thm:261112-intro} was with respect to the exact expression of
solutions of $P(\lm,\om_\a)$ with $\lm\leq \ov{\lm}$ and $\a$ small and not just with the estimates one
can get by using the sub-supersolutions just found in Theorem \ref{t1-intro}.\\

At this point (see section \ref{sec:entropy} for details) we can
conclude that indeed $S(E)\equiv \mathcal{S}(\de(\ull))\left.\right|_{\lm=\widehat{\lm}(E)}$ in
$\left[\frac{\a}{8\pi},\widehat{E}_\a\right]$.
In particular we conclude that $S(E)$ is also smooth
in $\left[\frac{\a}{8\pi},\widehat{E}_\a\right]$ and
by using the asymptotic expansions \rife{231112.3-intro}, \rife{231112.3.a-intro} and \rife{231112.3.b-intro}
and the above mentioned explicit expressions \rife{phi0-intro} and \rife{080213.1-intro}
we are eventually able to evaluate $\frac{d^2 S(E)}{d E^2}$ in the case $\om=\om_\a$ and
$E\in\left[\frac{\a}{8\pi},\widehat{E}_\a\right]$. Indeed, we have
\bpr\label{pr:entropy-intro}
Let $\om=\om_\a$, $E\in\left[\frac{\a}{8\pi},\widehat{E}_\a\right]$ and $\a<\tilde{\a}_1$ as defined in Theorem
\ref{unique:2-intro}. Then we have

$$
\frac{d^2 S(E)}{d E^2}=-11\left(\frac{48\pi^{2}}{\a^2}\right)+
\mbox{\rm O}\left(\frac{1}{\a}\right).
$$
\epr

In other words, we conclude that the branch of "small energy" solutions of $P(\lm,\om_\a)$ with $\lm\leq\ov{\lm}$
corresponds, for $\a$ small enough, to a range of energies where $S$ is concave.\\

\brm It can can be shown, of course with the necessary minor modifications,
that the monotonicity of the energy as a function of $\lm$ in Theorem \ref{unique:2-intro}, the
asymptotic expansion of the solution $\ull$ in Theorem \ref{thm:261112-intro} as well as
the concavity of the entropy in Proposition \ref{pr:entropy-intro}  still hold whenever
$\om$ is a regular domain such that
$\{\a^2x^2+y^2\leq \beta_-^2\}\subset\Omega\subset\{\a^2x^2+y^2\leq\beta_+^2\}$,
with $\frac{\bt^2_-}{\bt^2_+}=c$ and $c\in(\un{c},1]$ for some $\un{c}$ close enough to $1^{-}$.\\
The proofs of these results can be obtained with minor changes
by a step-by-step adaptation of those provided here. We will not discuss them
in particular because it seems that they do not provide any other useful
insight while they surely require a lot of additional technicalities.
\erm

\subsection{Open problems}

We conclude this introduction with a conjecture and an open problem.\\
It is well known that $S(E)$ is not concave (see NEQ-(iii) above) for $E>E_c$ and that solutions
of the MVP (see NEQ-(iii) and \rife{blow-up2} above) blow up as $E\to +\infty$.
Concerning this point we have the following:\\

{\bf Conjecture:} Let $\om$ be a convex domain of the second kind. There exists one and only one branch of solutions
$\ull$ which satisfies $\rife{blow-up2}$ and in particular there exists $E_{\om}>E_c$ such that $S(E)$ is convex in
$(E_\om,+\infty)$.\\

In particular uniqueness of blow-up solutions would imply that they coincide
(at least in a small right neighborhood of $8\pi$) with the set of
mountain-pass type solutions found in Theorem \ref{mp-intro}, see Remark \ref{metastable}.\\

Then we pose the following problem (see also fig.4 in \cite{clmp2}):\\

{\bf Open Problems:} Let us assume that either the above conjecture is true or that
$\om$ is a convex domain of the second kind for which we can find $E_{\om}>E_c$ such that $S(E)$ is convex in
$(E_\om,+\infty)$. Is it true that
the entropy has only one inflection point? If not, under which conditions (if any)
the entropy has only one inflection point?\\
In particular, is it true that the global branch of solutions of $P(\lm,\om_\a)$ with $\a$ small enough
has just one bending point, no bifurcation points and it is connected with the blow-up solution's branch
as $\lm\searrow (8\pi)^+$ (as for $Q(\mu,\om)$ on nearly circular domains \cite{suz0})? Can we answer this question at least on convex, regular and symmetric domains?\\

Of course, these properties do not hold on general simply connected domains. For example,
there should be no reason to expect the energy to be a generally injective function of $\lm$
(see for example fig.5 in \cite{clmp2}).
Moreover, some well known numerical results \cite{PW} suggest that bifurcation
points can exist on the bifurcation diagram of $P(\lm,\om)$ on (symmetric and/or non symmetric) non convex domains.
It seems however that the very rich structure of those bifurcation diagrams \cite{PW}
is inherited by solutions sharing either multiple peaks or just a single peak but which may
be located at different points. The typical example of such kind of blow-up behavior
is observed on dumbbell shaped domains, see for example \cite{EGP}.\\
On the other side, there are more lucky situations, such as on convex domains, where $k-$peaks
solutions with $k\geq 2$ do not exist (as shown in \cite{GrT}). Moreover it is well known (see for example \cite{Gu})
that if $\om$ is convex then the Robin function $H_\om(x,x)$ is strictly concave and thus admits
one and only one critical point, which of course coincides with the absolute maximum. This rules out the possibility
of having more than one single peak blow-up solution.\\
So far, it seems that in particular the global connectivity of the solution's
branch is known only for domains which are close in $C^2$-norm to a disk, see \cite{suz0}.\\
Of course, if (say in case $\om=\om_\a$ with $\a$ small enough) the entropy really has
just one inflection point, then it will coincide with the point on the continuation of $\mathcal{G}_{\a,1}$
where the first eigenvalue of the linearized problem for $P(\lm,\om_\a)$ will finally vanish.
However, in this situation we cannot use the standard results (see for example \cite{suzB}) which
in the classical cases show that this point must necessarily be a bending point. This is due to the peculiar form of the
linearized problem for $P(\lm,\om)$, see \rife{7.1} below, which implies for example that in general neither
the first eigenvalue can be assumed to be simple nor the first eigenfunction to be positive. This is not a mere
technical problem and indeed an explicit example of a sign changing first eigenfunction in a similar
situation can be found in Appendix D in \cite{B2}.\\
In any case we think that this topic deserves a separate discussion and that
it should be already very interesting to set up the problem on some
symmetric and convex domain of the second kind such as thin ellipses and/or rectangles.\\\\

This paper is organized as follows. In section \ref{sec:unique} we prove Theorem \ref{unique:1-intro}. In section
\ref{sec1} we prove Theorem \ref{t1-intro}. In section \ref{sec2} we prove Proposition \ref{pr3} by using a result
concerning the first eigenvalue of the linearization of $P(\lm,\om)$ around those solutions found in Theorem
\ref{t1-intro}, see Proposition \ref{pr2}. Section \ref{sec:mp} is devoted to the proof of Theorem \ref{mp-intro}.
Section \ref{s5} is concerned with the proofs of Theorems \ref{thm:261112-intro} and \ref{unique:2-intro}.
Finally section \ref{sec:entropy} is devoted to the proof of Proposition \ref{pr:entropy-intro}. Some technical
evaluations are left to the Appendix.

\bigskip
\bigskip

{\bf Acknowledgements.}\\
We wish to express our warmest thanks to Prof. M. Lassak for letting us know about His proof \cite{Lassek-priv}
of Theorem \ref{tlassek}. We are also indebted with Prof. G. Tarantello for suggesting the multiplicity
result Theorem \ref{mp-intro} and to Prof. A. Malchiodi for His suggestion about uniqueness of one-peak
blow-up solutions for the (MFE) and for His encouragement in our attempts to prove Theorems \ref{t1-intro}(b) and
\ref{unique:1-intro}.

\section{A uniqueness result for solutions of $P(\lm,\om)$.}\label{sec:unique}
The aim of this section is to obtain a uniqueness result for solutions of $P(\lm,\om)$ with
finite energy $\mathcal{E}(\de_{\scp \lm})\leq \ov{E}$ (see \rife{dedef}) on domains
chosen as in Theorem \ref{unique:1-intro}.

\bigskip

{\it The proof of Theorem \ref{unique:1-intro}}.\\
We will need an a priori estimate for solutions of $P(\lm,\om)$
with a uniform constant $\ov{C}$ which does not depend neither on $u$ nor on the domain $\om$. This is why we
do not follow the standard route which is
widely used (under some additional regularity assumption on $\pa\om$, see for example \cite{CCL})
in case the domain is fixed. In that case in fact one needs to prove that blow-up points
(in the sense of Brezis-Merle \cite{bm}) cannot converge to the boundary. A detailed discussion of
this point in our situation would be not only more tricky (since we do not fix $\om$)
but also really counterproductive, since
instead, by using the energy bound \rife{endef2-intro}, our argument
yields the needed estimate with the weakest possible regularity
assumptions about $\pa \om$ (i.e. $\om$ simple) see definition \ref{defsimp}.\\
The underlying idea is to use the dilation invariance (see Remark \ref{rem3.1})
of $P(\lm,\om)$ to show that even if a blow-up "bubble" converges to the boundary, then its energy must be unbounded.
More exactly we have:

\ble\label{lem:231112}
Let $\ov{\lm}\geq 8\pi$ and $\ov{E}\geq 1$ be fixed. There exists $\ov{C}=\ov{C}(\ov{\lm},\ov{E})$
such that for any simple domain $\om$ and for all solutions of
$P(\lm,\om)$ such that $\lm\leq \ov{\lm}$ and $\mathcal{E}(\de_{\scp \lm})\leq \ov{E}$ it holds
$\|u_{\scp \lm}\|_\i\leq \ov{C}$. In particular $\ov{C}$ does not depend neither on
$u$ nor on $\om$.\\
\ele
\proof In view of Remark \ref{solregdef} we can assume $u$ to be a classical solution
of $P(\lm,\om)$.\\
We argue by contradiction and suppose that there exists a sequence of simple domains $\{\om_n\}$ and a sequence of
positive numbers $\{\lm_n\}$ such that $\sup\limits_{\enne}\lm_n\leq\ov{\lm}$ and
there exists a sequence of solutions $\{u_n\}$ for $P(\lm_n,\om_n)$ such that
$$
\mathcal{E}(\de(u_n))\leq \ov{E},
$$
and there exists a sequence of points $\{x_n\}$ such that $x_n\in \om_n$ $\fo\,n\in\enne$ and
$$
u_n(x_n)=\max\limits_{\om_n}u_n\to+\infty.
$$
Of course, we have used here the fact that the maximum principle ensures that any solution for $P(\lm_n,\om_n)$
is nonnegative.\\
Since the problem is translation invariant we can assume without loss of generality that
$$
x_n\equiv 0,\;\fo n\in\enne.
$$

Let us set
$$
d_n:=\dist(0,\pa\om_n),
$$
and define
$$
w_{n,0}(y)=u_n\left(\frac{d_n}{2}y\right),\quad y\in \om_{n,0}:=\left\{y\in\R^2\,:\,\frac{d_n}{2}y\in\om_n\right\}.
$$

Clearly we have
\beq\label{100213.2}
B_1(0)\Subset \om_{n,0}
\eeq
and in particular (see Remark \ref{rem3.1}) $w_{n,0}$ is a solution of $P(\lm_n,\om_{n,0})$  which therefore satisfies
\beq\label{100213.3}
w_{n,0}(0)=u_n(0)=\max\limits_{\om_{n,0}}w_{n,0}\to+\infty.
\eeq

Let us set
$$
\mu_{n,0}:=\lm_n\left(\,\int\limits_{\om_{n,0}} e^{w_{n,0}}\right)^{-1}.
$$

We claim that:\\
{\bf Claim:} $w_{n,0}(0)+\log{\mu_{n,0}}\to+\infty$.\\
We argue by contradiction and observe that if the claim were false, then, in view of \rife{100213.3} we would have

$$
\graf{
-\Delta w_{n,0} \leq C_0  & \mbox{in}\quad \om_{n,0}\\
w_{n,0} =0 & \mbox{on}\quad \pa\om_{n,0}
}
$$
foe some $C_0>0$. For any $n\in\enne$ we can choose $R_n>0$ such that $\om_{n,0}\subset B_{R_n}$ and let
$$
\varphi_{n}(y)=\frac{C_0}{R_n^2}(R_n^2-|y|^2),\;y\in B_{R_n}
$$
be the unique solution of
$$
\graf{
-\Delta \varphi_n = C_0  & \mbox{in}\quad B_{R_n}\\
\varphi_n =0 & \mbox{on}\quad \pa B_{R_n}
}
$$

Clearly, by the maximum principle we have $w_{n,0}(0)\leq \varphi_n(0)=C_0$, which is a contradiction to
\rife{100213.3}.\fineproof

\bigskip

Therefore we see that the function $w_{n,1}(y)=w_{n,0}(y)+\log{\mu_{n,0}}$ satisfies

$$
\graf{
-\Delta w_{n,1} =e^{w_{n,1}}\quad \mbox{in}\quad B_1\\
\int\limits_{B_1}e^{w_{n,1}}\leq \ov{\lm}\\
w_{n,1}(0)=\max\limits_{B_1}w_{n,1}\to+\infty
}
$$

Hence we can apply the Brezis-Merle's result \cite{bm} as further improved by Li and Shafrir \cite{ls}
to conclude that there exists $r_0\in(0,1]$ such that
$$
e^{w_{n,1}}\rightharpoonup 8\pi m \delta_{p=0},\quad\mbox{in}\quad B_{2r_0},
$$
weakly in the sense of measures, where $m$ is a positive integer which satisfies
$1\leq m\leq\frac{\ov{\lm}}{8\pi}$. We remark that with a little extra work we could also
prove that the oscillation of $w_{n,1}$ is bounded on (say) $\pa B_{r_0}$ and hence in particular
obtain the desired contradiction by using the Li's result \cite{yy}. We will not purse this approach here
since we can come up with the desired conclusion just setting
\beq\label{110213.2}
\dt_{n,0}^2:=e^{-w_{n,1}(0)}\to 0,
\eeq
and use the by now standard blow-up argument in \cite{ls}. It shows that there exists a subsequence (which we will not relabel) such that
$$
w_n(z)=w_{n,1}(\dt_{n,0}z)-w_{n,1}(0),\quad |z|<(\dt_{n,0})^{-1},
$$
satisfies

\beq\label{110213.1}
w_n(z)\to w(z),\quad \mbox{in}\quad C^2_{\rm loc}(\R^2),
\eeq
where
\beq\label{110213.0}
w(z)=2\log{\frac{1}{(1+\frac{1}{8}|z|^2)}},\qquad \int\limits_{\R^2} e^w=8\pi.
\eeq

At this point we observe that, in view of the translation and dilation invariance of the energy we have
$$
\int\limits_{\om_{n,0}}\left| \nabla w_{n,1}\right|^2=\int\limits_{\om_{n,0}}\left| \nabla w_{n,0}\right|^2=
\int\limits_{\om_{n}}\left| \nabla u_{n}\right|^2=2\lm_n^2\mathcal{E}(\de(u_{n}))\leq 2\ov{\lm}^2\ov{E},
$$
so that, by using \rife{100213.3} and \rife{110213.2}, we should have,
$$
2\ov{\lm}^2\ov{E}\geq \int\limits_{\om_{n}}\left| \nabla u_{n}\right|^2=\lm_n\int\limits_{\om_{n}}\de(u_{n})u_{n}=
\lm_n\int\limits_{\om_{n,0}}\de(w_{n,0})w_{n,0}>\lm_n\int\limits_{B_{R\dt_{n,0}}}\de(w_{n,0})w_{n,0}=
$$
$$
\int\limits_{B_{R\dt_{n,0}}}e^{w_{n,1}}(w_{n,1}-\log\mu_{n,0})=\int\limits_{B_R}e^{w_{n}}(w_{n}+w_{n,1}(0)-\log\mu_{n,0})=
$$
$$
\int\limits_{B_R}e^{w_{n}}w_{n}+u_{n}(0)\int\limits_{B_R}e^{w_{n}},
$$
for any $R\geq 1$ and for any $n\in\enne$,  which is clearly in contradiction with \rife{100213.3} and \rife{110213.1}, \rife{110213.0}.
We refer to Lemma 3.1 in \cite{bl} for a proof of the fact that the Gauss-Green formula
$\int\limits_{\om_{n}}\left| \nabla u_{n}\right|^2=\lm_n\int\limits_{\om_{n}}\de(u_{n})u_{n}$ holds
on domains which are only assumed to be simple.
\fineproof

\bigskip

{\it The proof of Theorem \ref{unique:1-intro} completed}.\\
We first prove part (b).\\
We argue by contradiction and suppose that there exists a sequence of open, bounded and convex domains
$\{\om_{n,0}\}$ such that
\beq\label{hpdiv}
N(\om_{n,0})=\frac{L^2(\pa \om_{n,0})}{A(\om_{n,0})}> n,
\eeq
and a sequence of positive numbers $\{\lm_n\}$ such that $\sup\limits_{\enne}\lm_n\leq\ov{\lm}$,
such that for any $n\in\enne$ there exist at least two solutions $u_{n,1}$ and $u_{n,2}$
for $P(\lm_n,\om_{n,0})$ such that
\beq\label{110213.5}
\mathcal{E}(\de(u_{n,i}))\leq \ov{E},\quad i=1,2.
\eeq

In view of Theorems \ref{tjohn} and \ref{tlassek} we see that for each $n\in\enne$ there exist two concentric and
omotetic ellipses such that
\beq\label{Lassek.0}
\mathbb{E}_{n,-}\subseteq\om_{n,0}\subseteq\mathbb{E}_{n,+}
\eeq
and
\beq\label{Lassek}
\frac{A(\mathbb{E}_{n,+})}{A(\mathbb{E}_{n,-})}=4.
\eeq

Since $P(\lm,\om)$ and \rife{110213.5} are both rotational, translational and dilation invariant, then,
in view of Remark \ref{rem3.1}, we can assume
without loss of generality that for each $n\in\enne$
\beq\label{110213.10}
\mathbb{E}_{n,+}=\om_{\a_n},\quad\mbox{for some}\quad \a_n>0.
\eeq

Clearly we have

$$
N(\mathbb{E}_{n,+})=\frac{L^2(\pa \mathbb{E}_{n,+})}{A(\mathbb{E}_{n,+})}=
\frac{\L^2(\pa \mathbb{E}_{n,+})}{4A(\mathbb{E}_{n,-})}\geq
\frac{L^2(\pa \mathbb{E}_{n,+})}{4A(\om_{n,0})}\geq \frac14 N(\om_{n,0})>\frac{n}{4}.
$$

Therefore, since in view of \rife{110213.10} we have $L^2(\pa \mathbb{E}_{n,+})\leq \frac{4\pi^2}{\a_n^2}$ and
$A(\mathbb{E}_{n,+})=\frac{\pi}{\a_n}$, then we also conclude that
$$
\frac{n}{4}<N(\mathbb{E}_{n,+})\leq \frac{4\pi^2}{\a_n^2}\frac{\a_n}{\pi},
$$
that is
\beq\label{110213.11}
\a_n<\frac{16\pi}{n}\,.
\eeq

At this point we observe that
\beq\label{110213.21}
\sg_{n,0}:=\inf\left\{
\left.\frac{\int\limits_{\om_{n,0}} \left|\nabla \vp\right|^2\,dx}
{\int\limits_{\om_{n,0}}\vp^2\,dx}\;\right|
\,\vp\in H^{1}_0(\om_{n,0})\right\}\geq  2(1+\a_n)>2,
\eeq
which easily follows from the fact that $\sg_{n,0}\geq \sg_n$, where
$$
\sg_n:=\inf\left\{
\left.\frac{\int\limits_{\om_{\a_n}} \left|\nabla \vp\right|^2\,dx}
{\int\limits_{\om_{\a_n}}\vp^2\,dx}\;\right|
\,\vp\in H^{1}_0(\om_{\a_n})\right\},
$$
see \rife{140213.0} and \rife{14} below for further details.\\
Hence, by using \rife{110213.21}, we conclude that
$$
2\int\limits_{\om_{n,0}}\left|u_{n,1}-u_{n,2}\right|^2\leq
\int\limits_{\om_{n,0}}\left|\nabla(u_{n,1}-u_{n,2})\right|^2=
\lm_n\int\limits_{\om_{n,0}}(\de(u_{n,1})-\de(u_{n,2}))(u_{n,1}-u_{n,2}).
$$
Let us write
$$
\int\limits_{\om_{n,0}}(\de(u_{n,1})-\de(u_{n,2}))(u_{n,1}-u_{n,2})=I_{1,n}+I_{2,n},
$$
where
$$
I_{1,n}=
\int\limits_{\om_{n,0}}\frac{e^{u_{n,1}}-e^{u_{n,2}}}{\int\limits_{\om_{n,0}} e^{u_{n,1}}}(u_{n,1}-u_{n,2}),\quad
I_{2,n}=
\int\limits_{\om_{n,0}}e^{u_{n,2}}
\left(\frac{1}{\int\limits_{\om_{n,0}} e^{u_{n,1}}}-\frac{1}{\int\limits_{\om_{n,0}} e^{u_{n,2}}}\right)(u_{n,1}-u_{n,2}).
$$
It follows from Lemma \ref{lem:231112} (which of course can be applied since any open, bounded and convex
domain is simple according to Definition \ref{defsimp}) and the fact that solutions of $P(\lm,\om)$
are non negative that, by using also \rife{Lassek}, we can estimate these two integrals as follows
$$
\left|I_{1,n}\right|\leq \int\limits_{\om_{n,0}}\frac{e^{\ov{u}_n}}{A(\om_{n,0})}|u_{n,1}-u_{n,2}|^2\leq
\int\limits_{\om_\a}\frac{e^{\ov{C}}}{A(\mathbb{E}_{n,-})}|u_{n,1}-u_{n,2}|^2=
\frac{4e^{\ov{C}}}{\pi}\a_n \int\limits_{\om_{n,0}}|u_{n,1}-u_{n,2}|^2,
$$
and similarly,
$$
\left|I_{2,n}\right|\leq
\int\limits_{\om_{n,0}}e^{u_{n,2}}|u_{n,1}-u_{n,2}|\left(\int\limits_{\om_{n,0}}\frac{e^{\ov{u}_n}}
{\left(\int\limits_{\om_{n,0}}e^{\ov{u}_n}\right)^2}(u_{n,1}-u_{n,2})\right)\leq
\frac{e^{2\ov{C}}}{A^2(\om_{n,0})}\left(\int\limits_{\om_{n,0}}|u_{n,1}-u_{n,2}|\right)^2\leq
$$

$$
\frac{4 e^{2\ov{C}}}{\pi A(\om_{n,0})}\a_n\left(\int\limits_{\om_{n,0}}|u_{n,1}-u_{n,2}|\right)^2\leq
\frac{4 e^{2\ov{C}}}{\pi}\a_n\int\limits_{\om_{n,0}}|u_{n,1}-u_{n,2}|^2,
$$
where $\ov{u}_n$ is a suitable function which satisfies
$\ov{u}_n\in(\min\{u_{n,1},u_{n,2}\}, \max\{u_{n,1},u_{n,2}\})$.\\
Plugging these estimates together we conclude that
$$
\int\limits_{\om_{n,0}}|u_{n,1}-u_{n,2}|^2\leq \lm_n\a_n\frac{8e^{2\ov{C}}}{\pi}\int\limits_{\om_{n,0}}|u_{n,1}-u_{n,2}|^2,
$$

which is of course a contradiction to \rife{110213.11}.
This contradiction shows that in fact there exists at most one solution under the given assumptions and
concludes the proof of part (b) of the statement.\\

As for part (a) it is easy to adapt the argument by contradiction used above just by
replacing the assumption of divergent isoperimetric ratio in \rife{hpdiv} with that of
the existence of $\a_n\searrow 0^+$ and $0<\beta_{-,n}\leq \beta_{+,n}<+\i$ such that
$$
\mathbb{E}_{n,-}:=\{\a_n^2x^2+y^2\leq \beta_{-,n}^2\}\subseteq\om_{n,0}\subseteq
\{\a_n^2x^2+y^2\leq \beta_{+,n}^2\}=:\mathbb{E}_{n,+},\quad \frac{\beta_{-,n}}{\beta_{+,n}}=c,\quad \fo\,n\in\enne.
$$

In particular we see that this time we already have (by assumption)
the needed concentric omotetic ellipses (as in \rife{Lassek.0}) which in this case satisfy
$$
\frac{A(\mathbb{E}_{n,+})}{A(\mathbb{E}_{n,-})}=\frac{\beta_{+,n}^2}{\beta_{-,n}^2}=c^2.
$$

At this point, since of course Lemma \ref{lem:231112} can be applied to the situation at hand, the proof
can be worked out as above with minor changes.\fineproof

\section{Solutions of supercritical Mean Field Equations on thin domains}\label{sec1}

In this section we prove Theorem \ref{t1-intro}. Indeed, we will
construct a branch of solutions of $P(\lm,\om_\a)$ which for $\a$ small enough extends up to some value
${\lm}_{\a}\geq \frac{4\pi}{7\a}$, and more generally we obtain the same statement on any domain
$\Omega$ lying between two concentric and similar \virg{thin} ellipses.
Thus, in particular we recover the result for convex domains having a large isoperimetric ratio.
To achieve our goal, we consider the auxiliary problem $Q(\mu,\om)$ (see \S \ref{ss1.1})
and make use of a well known result \cite{clsw} whose statement calls up for the following:
\bdf\label{defss}
A function $u$ is said to be a \un{subsolution}(\un{supersolution}) of $Q(\mu,\om)$
if $u\in C^0(\ov{\om})$ and
\beq\label{sssyst}
\graf{
\ino(-\Delta\varphi) u\leq(\geq) \mu \,{\dsp e^u}\varphi & \mbox{in}\quad \om\\\\
u\leq(\geq)0 & \mbox{on}\quad \pa\om
},\qquad\fo\varphi\in C^{\infty}_0(\om),\varphi\geq 0.
\eeq
\edf

\bigskip

\bte[Sub-Supersolutions method, \cite{clsw}]\label{t0}
Let $\om$ be simple.
Suppose that, for fixed $\mu>0$, there exist
a subsolution $\subu_\mu$ and a supersolution $\supu_\mu$ of $Q(\mu,\om)$.
If $\subu_\mu\leq\supu_\mu$ in $\om$, then
$Q(\mu,\om)$ admits a classical solution $u=u_\mu\in C^{2}(\om)\cap C^0(\ov{\om})$
which moreover satisfies $\subu_\mu\leq u_\mu \leq \supu_\mu$.
\ete

\proof
We use the existence Theorem in \cite{clsw}, where
the domain $\om$ is just assumed to be regular with respect to the Laplacian (see \cite{GT}, p. 25). It is well
known that any simple domain satisfies this assumption (see \cite{GT}, p. 26).
Therefore we can apply the result in \cite{clsw} which yields the existence of a function
$u_\mu\in C^0(\ov{\om})$ which satisfies $\subu_\mu\leq u_\mu \leq \supu_\mu$ and moreover satisfies
\rife{sssyst} for \un{all} $\varphi\in C^{\infty}_0(\om)$ with the equality sign replacing the corresponding
inequalities. Hence in particular $u_\mu$
is a distributional solution of the equation in $Q(\mu,\om)$. Therefore the Brezis-Merle \cite{bm} theory of
distributional solutions of Liouville type equations shows that it is also locally bounded and then
standard elliptic regularity theory shows that $u_\mu\in C^{2}(\om)$ is a classical solution of $Q(\mu,\om)$ as well.
We insist about the fact that the continuity up to the boundary is a byproduct of the result in \cite{clsw},
which indeed yields a distributional solution $u_\mu\in C^0(\ov{\om})$.
\fineproof

\proof[Proof of Theorem \ref{t1-intro}(a)]
For fixed $c\in(0,1]$ and in view of Remark \ref{rem3.1}
we can assume without loss of generality that
$$
\Omega_{\rho,c}:=\{\rho^2 x^2+y^2\leq c\}\subseteq\Omega\subseteq\{\rho^2 x^2+y^2\leq 1\}=:\Omega_\rho.
$$

Let us define
\beq\label{2}
v_{\a,\ga}=2\log{\left(\frac{1+\ga^2}{1+\ga^2(\a^2x^2+y^2)}\right)},\quad (x,y)\in \om_\a.
\eeq

A straightforward evaluation shows that $v_{\a,\ga}$ satisfies
\beq\label{3}
\graf{
-\Delta v_{\a,\ga}=V_{\a,\ga}{\dsp e^{v_{\a,\ga}}} & \mbox{in}\quad \om_\a\\
v_{\a,\ga}=0 & \mbox{on}\quad \pa\om_\a,
}
\eeq
where
\beq\label{4}
V_{\a,\ga}(x,y)=\frac{4\ga^2}{(1+\ga^2)^2}\left(1+\a^2+\ga^2(1-\a^2)(\a^2x^2-y^2)\right)
\eeq

Since
$$
V_{\a,\ga}(x,y)\geq g_{+}(\ga,\a):=\frac{4\ga^2}{(1+\ga^2)^2}\left(1+\a^2+\ga^2(\a^2-1)\right),\quad \forall (x,y)\in \om_\a,
$$
we easily verify that $v_{\a,\ga}$ is a classical supersolution and in particular a
supersolution (according to the above definition) of $Q(\mu,\om)$ whenever
\beq\label{5}
\mu\leq g_{+}(\ga,\a).
\eeq
For fixed $\a\in(0,1)$, the function $h_\a(t)=g_{+}(\sqrt{t},\a)$ satisfies
$h_\a(0)=0=h_\a\left(\frac{1+\a^2}{1-\a^2}\right)$, is strictly increasing in $\left(0,\frac{1+\a^2}{3-\a^2}\right)$ and
strictly decreasing in $\left(\frac{1+\a^2}{3-\a^2},\frac{1+\a^2}{1-\a^2}\right)$.
Therefore, putting $\overline{\ga}_\a^2=\frac{1+\a^2}{3-\a^2}$ and
$\overline{\mu}_\a:=h_\a\left(\overline{\ga}_\a^2\right)\equiv g_{+}(\overline{\ga}_\a,\a)\equiv \frac{(\a^2+1)^2}{2}$,
we see in particular that for each $\mu\in(0,\overline{\mu}_\a]$ there exists a unique
$\ga^{+}_\a\in\left(0,\overline{\ga}_\a\right]$ such that $g_{+}(\ga^{+}_\a,\a)=\mu$ and $v_{\a,\ga^{+}_\a}$ is
a supersolution of $Q(\mu,\om)$. Indeed we have
$$
\left(\ga^{+}_\a\right)^2=\left(\ga^{+}_\a(\mu)\right)^2=\frac{2(1+\a^2)-\mu-2\sqrt{(1+\a^2)^2-2\mu}}{\mu+4(1-\a^2)}.
$$

On the other hand let us consider
\beq\label{vmeno}
v_{\a,\ga,c}=\left\{
               \begin{array}{ll}
                 2\log{\left(\frac{1+\ga^2}{1+\tfrac{\ga^2}{c}(\a^2x^2+y^2)}\right)},& \hbox{$(x,y)\in\Omega_{\rho,c}$} \\
                 0, & \hbox{$(x,y)\in\Omega\setminus\Omega_{\rho,c}$.}
               \end{array}
             \right.
\eeq

Again a straightforward computation shows that $v_{\a,\ga,c}$ satisfies
$$
\left\{
  \begin{array}{ll}
-\Delta v_{\a,\ga,c}=V_{\a,\ga,c}{\dsp e^{v_{\a,\ga,c}}} & \mbox{in}\quad \Omega_{\a,c}\\
v_{\a,\ga,c}=0 & \mbox{on}\quad \pa\Omega_{\a,c},  \end{array}
\right.
$$
where
$$
V_{\a,\ga,c}(x,y)=\left\{
               \begin{array}{ll}
            \frac{4\ga^2}{c(1+\ga^2)^2}\left(1+\a^2+\frac{\ga^2}{c}(1-\a^2)(\a^2x^2-y^2)\right) & \mbox{in }\Omega_{\a,c} \\
                 0 & \mbox{in }\Omega\setminus\Omega_{\a,c}.
               \end{array}
             \right.
$$
Since
$$
V_{\a,\ga,c}(x,y)\leq g_{-}(\ga,\a,c):=\frac{4\ga^2}{c(1+\ga^2)^2}\left(1+\a^2+\ga^2(1-\a^2)\right),\quad \forall (x,y)\in \om,
$$
it is not difficult to check that $v_{\a,\ga,c}$ is a subsolution (according to the above definition) of
$Q(\mu,\om)$ whenever
\beq\label{6}
\mu\geq g_{-}(\ga,\a,c).
\eeq
For fixed $\a\in(0,1)$, the function $f_{\a,c}(t)=g_{-}(\sqrt{t},\a,c)$, $t\in(0,\overline{\ga}_\a^2]$ is strictly
increasing and
satisfies $f_{\a,c}(t)>h_\a(t)$.
Therefore, for each $\mu\in(0,\overline{\mu}_\a]$ there exists a unique
$\ga^{-}_{\a,c}\in\left(0,\overline{\ga}_\a\right)$ such that $g_{-}(\ga^{-}_{\a,c},\a,c)=\mu$,
$\ga^{-}_{\a,c}<\ga^{+}_\a$ and $v_{\a,\ga^{-}_{\a,c},c}$ is
a  subsolution of $Q(\mu,\om)$. Indeed we have
$$
\left(\ga^{-}_{\a,c}\right)^2=\left(\ga^{-}_{\a,c}(\mu)\right)^2=\frac{\mu c-2(1+\a^2)+2\sqrt{(1+\a^2)^2-2\a^2\mu c}}{4(1-\a^2)-\mu c}.
$$
In conclusion, since $\ga^{-}_{\a,c}(\mu)\leq \ga^{+}_\a(\mu)$ implies $v_{\a,\ga^{-}_{\a,c},c}\leq v_{\a, \ga^{+}_\a}$,
for fixed $\a\in (0,1)$ and for each $\mu\in(0,\overline{\mu}_\a]$ we can set
$$
\subu_{\mu}=v_{\a,\ga^{-}_{\a,c}(\mu),c},\quad \supu_\mu=v_{\a,\ga^{+}_\a(\mu)},
$$
to obtain (through Theorem \ref{t0}) a solution $u_{\a,\mu,c}$ for $Q(\mu,\om)$ which satisfies
\beq\label{susu}
v_{\a,\ga^{-}_{\a,c}(\mu),c}\leq u_{\a,\mu,c}\leq v_{\a,\ga^{+}_\a(\mu)},\quad \forall (x,y)\in\om.
\eeq
Any such a solution $u_{\a,\mu,c}$ therefore solves $P(\lm,\om)$ with $\lm=\lm_{\a,c}(\mu)$ satisfying
\bel\label{stimalambdamu1}
\lm=\lm_{\a,c}(\mu)=\mu\int\limits_{\om}e^{u_{\a,\mu,c}}\geq \mu\int\limits_{\om_{\a,c}}e^{v_{\a,\ga^{-}_{\a,c}(\mu),c}}=
\mu c\frac{\pi}{\a}\,(1+(\ga^{-}_{\a,c}(\mu))^2),
\eel
and
\bel\label{stimalambdamu2}
\lm=\lm_{\a,c}(\mu)=\mu\int\limits_{\om}e^{u_{\a,\mu,c}}\leq \mu\int\limits_{\om_\a}e^{v_{\a,\ga^{+}_\a(\mu)}}=
\mu\frac{\pi}{\a}\,(1+(\ga^{+}_\a(\mu))^2).
\eel

In the particular case $\mu=\overline{\mu}_\a$ we have $(\ga^{-}_{\a,c}(\overline{\mu}_\a))^2\equiv
\underline{\ga}_{\a,c}^2=(1+\a^2)
\frac{c-4+c\a^2+4\sqrt{1-c\a^2}}{8(1-\a^2)-c(1+\a^2)^2}$,
$\underline{\ga}_{\a,c}^2<\overline{\ga}_\a^2$, $(\ga^{+}_\a(\overline{\mu}_\a))^2\equiv\overline{\ga}_\a^2=(1+\a^2)
\frac{3-\a^2}{8(1-\a^2)+(1+\a^2)^2}$ and
$u_{\a,\overline{\mu}_\a,c}$ is a solution for $P(\lm_{\a,c}(\overline{\mu}_\a),\om)$, where
\bel\label{lambdasotto}
\lm_{\a,c}:=\lm_{\a,c}(\overline{\mu}_\a)\geq \underline{\lm}_{\a,c}=\frac{c(1+\a^2)^2}{2}\frac{\pi}{\a}(1+\underline{\ga}_{\a,c}^2)\simeq
\frac{4\pi c}{(8-c)\a},
\eel
and
\bel\label{lambdasopra}
\lm_{\a,c}:=\lm_{\a,c}(\overline{\mu}_\a)\leq \overline{\lm}_{\a}=\frac{(1+\a^2)^2}{2}\frac{\pi}{\a}
(1+\overline{\ga}_\a^2)\simeq\frac{11\pi}{16\a}
\eel
as $\a\to 0^+$. Moreover it is easy to verify that $\un{\lm}_{\a,c}$ is strictly decreasing at least for
for $\a\in(0,\frac 1{2\sqrt{10}}]$ and that there exists
$\underline{\a}_*(c)<\frac 1{2\sqrt{10}}$ such that
$\underline{\lm}_{\a,c}\geq 8\pi$ for any $\a\in(0,\underline{\a}_*(c)]$.
We also see that $\ov{\lm}_\a\to 4\pi^{-}$ as $\a\to 1^{-}$,
is strictly decreasing for $\a\in(0,\a_p]$ and strictly increasing
for $\a\in[\a_p,1)$ for some $\a_p\simeq 0.5$ and then it is straightforward to check that
there exists $\overline{\a}_*>\underline{\a}_*(c)$ satisfying $0.0702<\overline{\a}_*<0.0703$ such that
$\overline{\lm}_\a\geq 8\pi$ for any $\a\in(0,\overline{\a}_*]$.

Finally, since $\lm_{\a,c}(\mu)$ is continuous in $\mu$ and by using \eqref{stimalambdamu1} and \eqref{stimalambdamu2}
$$
0<\lm_{\a,c}(\mu)\leq\mu\frac{\pi}{\a}\,(1+(\ga^{+}_\a(\mu))^2)\stackrel{\textnormal{as $\mu\to0$}}{\longrightarrow} 0,
$$
we obtain the existence of a solution for $P(\lambda,\Omega)$ not only for
$\lambda=\lm_{\a,c}$, but for any $\lambda\in(0,\lm_{\a,c}]$ as well.\fineproof

\bigskip

\proof[Proof of Theorem \ref{t1-intro}(b)] If ${\bar N}$ exists, then Remark \ref{susp} shows that
it is strictly greater than $4\pi$.\\
In view of Remark \ref{rem3.1} we can assume without loss of generality that $L(\pa\Omega)=1$.
Let $E_1$ be the John maximal ellipse of $\Omega$,
then by Theorem \ref{tjohn} $E_2:=\{c_0+2(x-c_0):x\in E_1\}$, where $c_0$ is the center of $E_1$, contains $\Omega$.
Again by using Remark \ref{rem3.1} we can also assume that $c_0=0$ and in particular that $E_1$ and $E_2$ have the following form
$$
E_1=\left\{\frac{x^2}{a^2}+\frac{y^2}{b^2}=1\right\},\quad E_2=\left\{\frac{x^2}{a^2}+\frac{y^2}{b^2}=4\right\},
$$
where clearly we can suppose that $0<b\leq a$.

By virtue of Ramanujan's estimate of the perimeter of the ellipse \cite{Villarino}, namely:
$$
L(\partial E_1)\geq \pi\{(a+b)+\frac{3(a-b)^2}{10(a+b)+\sqrt{a^2+14ab+b^2}}\},
$$
being $E_1\subset\Omega$, $\Omega$ convex, and since $N(\om)=\frac{L^2(\pa\Omega)}{A(\Omega)}$, we derive the following inequalities:
\bel\label{stimeab1}
1=L(\pa\Omega)\geq L(\pa E_1)\geq(a+b)\pi;\quad\qquad \frac{1}{N(\om)}=\frac{A(\Omega)}{L^2(\pa\Omega)}=A(\Omega)\geq A(E_1)=\pi a b.
\eel
Moreover since $\Omega\subset E_2\subset R_{a,b}:=\{(x,y)\in\R^2\,|\, |x|\leq 2a,\,|y|\leq 2b\}$ we get
\bel\label{stimeab2}
1=L(\pa\Omega)\leq L(\pa E_2)\leq L(R_{a,b})=8(a+b),
\eel
and by using Theorem \ref{tlassek}
\bel\label{stimeab3}
\frac{1}{N(\om)}=A(\Omega)\leq \frac{3\sqrt3}{\pi}A(E_1)=3\sqrt3 ab.
\eel
To simplify the notation we set $N=N(\Omega)$. Collecting \eqref{stimeab1}, \eqref{stimeab2} and \eqref{stimeab3} we have
\bel\label{stimeab4}
\left\{
  \begin{array}{l}
   \frac1{3\sqrt3 N}\leq ab\leq \frac1{\pi N}\\
-b+\frac1{8}\leq a\leq -b+\frac1\pi,
  \end{array}
\right.
\eel
which in turn implies
$$
\left\{
  \begin{array}{l}
b^2-\frac b\pi+\frac1{3\sqrt3 N}\leq 0\\
b^2-\frac b8+\frac1{\pi N}\geq 0.
  \end{array}
\right.
$$
It is worth to notice that, since $a\geq b$ and $ab\leq\frac{1}{\pi N}$, if $N>\frac{64}{\pi}$ then $b<\frac18$.
Therefore solving the above system of inequalities, with $N>\frac{64}\pi$, we get
$$
\frac{1-\sqrt{1-\frac{4\pi^2}{3\sqrt3 N}}}{2\pi}\leq b\leq \frac{1-\sqrt{1-\frac{256}{\pi N}}}{16}.
$$
Next, for $N>\frac{512}{\pi}$, considering the Taylor formula of the square root and estimating
the second order reminder we derive
\bel\label{b}
\frac{\pi}{3\sqrt3 N}\leq\frac{\frac{2\pi^2}{3\sqrt3 N}+\frac18(\frac{4\pi^2}{3\sqrt3 N})^2}{2\pi}\leq b\leq
\frac{\frac{128}{\pi N}+\frac{1}{2\sqrt2}(\frac{256}{\pi N})^2}{16}=\frac{8}{\pi N}+\frac{1024\sqrt2}{\pi^2 N^2},
\eel
thus
\bel\label{a}
\frac18-\frac{8}{\pi N}-\frac{1024\sqrt2}{\pi^2 N^2}\leq a\leq\frac1\pi-\frac\pi{3\sqrt3 N}.
\eel
Combining \eqref{b} and \eqref{a}, we have
$$
\psi(N):=\frac{\pi^2}{3\sqrt3 N-\pi^2}\leq\frac {b}{a}\leq
\frac{64+\frac{8192\sqrt2}{\pi N}}{\pi N-64-\frac{8192\sqrt2}{\pi N}}=:\varphi(N).
$$
By definition of $E_1$ and $E_2$ we are in position to apply point (a) of this theorem with $c=\frac14$.
Let us fix $\bar N$ such that
$\frac{64+\frac{8192\sqrt2}{\pi \bar N}}{\pi \bar N-64-\frac{8192\sqrt2}{\pi \bar N}}=
\underline{\a}_*(\frac14)$. We point out that since $\underline{\a}_*(\frac14)\simeq0,0161$,
$\bar N>\frac{512}{\pi}$.

Then, for any $N\geq\bar N$, $\a_{\scp N}:=\frac{b}{a}\leq \underline{\a}_*(\frac{1}{4})$ and so we get
the existence of a solution $\ul$ to $P(\lm,\om)$ for any $\lm\leq\lm_{\textnormal{\tiny{$N$}}}$ where
$$
\underline{\Lambda}_\textnormal{\tiny{$N$}}:=\underline{\lambda}_\textnormal{\tiny{$\varphi(N),\frac14$}}
\leq\underline{\lambda}_\textnormal{\tiny{$\a_N,\frac14$}}
<\lm_\textnormal{\tiny{$N$}}<\overline{\lambda}_\textnormal{\tiny{$\a_N$}}
\leq\overline{\lambda}_\textnormal{\tiny{$\psi(N)$}}=:\overline{\Lambda}_\textnormal{\tiny{$N$}}.
$$

At last from \eqref{lambdasotto} and \eqref{lambdasopra} we obtain the desired estimates on
$\underline{\Lambda}_\textnormal{\tiny{$N$}}$ and $\overline{\Lambda}_\textnormal{\tiny{$N$}}$:
$$
\underline{\Lambda}_\textnormal{\tiny{$N$}}\simeq\frac{\pi^2N}{496}+O(1)
\qquad \overline{\Lambda}_\textnormal{\tiny{$N$}}\simeq\frac{33\sqrt3N}{16\pi}
+O(1)\quad \textnormal{as $N\to+\infty$.}
$$

\fineproof

\section{The eigenvalue problem}\label{sec2}

The aim of this section is to prove Proposition \ref{pr2} below which yields positivity of the
first eigenvalue for the linearization of $P(\lm,\om)$. Among other things,
with the aid of Proposition \ref{pr2} we have:

\proof[Proof of Proposition \ref{pr3}]
Let $\mathcal{G}_{\a,c},\mathcal{G}_{N}$ denote the set of pairs of parameter-solutions
for $P(\lm,\om)$ found in Theorem \ref{t1-intro}. Since the linearized problem for
$P(\lm,\om)$ corresponds to the kernel equation for the second variation of $J_\lm$,
then the conclusions of Proposition \ref{pr3}
are an immediate consequence of Proposition \ref{pr2} below and the uniqueness results in \cite{CCL}.
\fineproof

\bigskip

Putting
$$
\de=\de(u)=\dfrac{e^u}{\int\limits_{\om} e^u},\quad \mbox{and}\quad <f>_\de=\int\limits_{\om}\de(u)f,
$$
then the linearized problem for $P(\lm,\om)$ takes the form
\bel\label{7.1}
\left\{
  \begin{array}{ll}
     -\Delta \vp - \lm\de(u)\vp +\lm \de(u)<\vp>_\de =0 & \mbox{in}\ \ \om \\
    \vp =0 & \mbox{on}\ \ \pa\om.
  \end{array}
\right.
\eel

\bpr\label{pr2}
For fixed $c\in(0,1]$, let $\om$ be a regular domain
such that $\{\a^2x^2+y^2\leq \beta_-^2\}\subset\Omega\subset\{\a^2x^2+y^2\leq\beta_+^2\}$,
with $\frac{\bt^2_-}{\bt^2_+}=c$. For any $\a\in(0,\underline{\a}_*(c)]$
let $u=\ul\equiv u_{\a,\mu,c}$ be a solution of $P(\lm,\om)$ and of $Q(\mu,\om)$
with $\lm=\mu\int\limits_{\om}e^u$ as obtained in Theorem \ref{t1-intro}(a) for $\lambda\in[0,\lm_{\a,c}]$.
Then \eqref{7.1} has only the trivial solution and in particular
the first eigenvalues of the
linearized problems for $P(\lm,\om)$ and $Q(\mu,\om)$ at $u=\ul\equiv u_{\a,\mu,c}$ respectively are strictly positive.\\
Moreover, let $\om$ be a regular and convex domain with $N(\om)>\bar N$ as defined in Theorem \ref{t1-intro}(b)
and $\ul$ be a solution of $P(\lm,\om)$ and of $Q(\mu,\om)$
for $0\leq \lm=\mu\int\limits_{\om}e^{\ul}\leq \lm_{\textnormal{\tiny{$N$}}}$
as obtained therein. Then the first eigenvalues of the
linearized problems for $P(\lm,\om)$ and $Q(\mu,\om)$ at $u=\ul$ are strictly positive.
\epr

\brm\label{branch}
As far as one is concerned with problem $Q(\mu,\om)$, then it is well known (see for example \cite{suzB}) that
is well defined (and unique) the extremal (classical) solution $v_*$ which corresponds to the extremal value $\mu_*$ such that no solutions exists
for $\mu>\mu_*$ and the bifurcation diagram has a bending point at $(\mu_*,v_*)$. In particular the first eigenvalue
of the linearized problem for $Q(\mu,\om)$ is zero at $\mu_*$.\\
The reasons why we have strictly positive first eigenvalues for $\lm\leq\lm_{\a,c}$ are:\\
(-) as it will be shown in the proof below, the first eigenvalue of the linearized problem for $P(\lm,\om)$ (say $\tau_1$)
is always greater or equal to the first eigenvalue (which we will denote by $\nu_0$)
of the linearized problem for $Q(\mu,\om)$ and we will use the latter to estimate both;\\
(-) the value of $\mu$ corresponding to $\lm_{\a,c}$,
which is defined implicitly via $\lm=\mu\int\limits_{\om_\a}e^{u_{\scp \mu}}$, is less than $\mu_*$.\\
\erm

\proof
We will use the fact that (see \cite{bm}, \cite{CCL} and Remark \ref{solregdef} above)
if $u$ solves $P(\lm,\om)$ then there exists $C=C(\om,\lm,u)>0$ such that
$$
\frac{1}{C}\leq \de(u)\leq C.
$$

Letting $H\equiv H^{1}_0(\om)$ and
$$
\mathcal{L}(\phi,\psi)= \int\limits_{\om} \left(\nabla \phi\cdot\nabla \psi\right)-
\lm\int\limits_{\om} \de(u)\phi\psi +
\lm\left(\int\limits_{\om} \de(u)\phi\right)\left(\int\limits_{\om} \de(u)\psi\right),\;(\phi,\psi)\in H\times H,
$$
then by definition $\vp\in H$ is a weak solution of \eqref{7.1} if
$$
\mathcal{L}(\vp,\psi)=0,\quad \forall\, \psi \in H.
$$
We define $\tau\in\erre$ to be an eigenvalue of the operator
$$
L[\vp]:=-\Delta\vp  - \lm\de(u) (\vp -< \vp >_\de),\quad \vp\in H,
$$
if there exists a weak solution $\phi_0\in H\setminus \{0\}$ of the linear problem

\beq\label{7.2}
 -\Delta \phi_0 - \lm\de(u)\phi_0+\lm \de(u)<\phi_0>_\de =\tau \de(u)\phi_0\ \ \mbox{in}\ \ \om,
\eeq
that is, if
$$
\mathcal{L}(\phi_0,\psi)=\tau\int\limits_{\om} \de(u)\phi_0\psi,\quad \forall\, \psi \in H.
$$

Standard arguments show that the eigenvalues form an unbounded (from above) sequence
$$
\tau_1\leq \tau_2\leq\cdots\leq \tau_n\cdots,
$$

with finite dimensional eigenspaces  (although the first eigenfunction cannot be assumed to be neither positive nor
simple in this situation).

Let us define
$$
Q(\phi)=\frac{\mathcal{L}(\phi,\phi)}{<\phi^2>_\de}=
\frac{\int\limits_{\om} \left|\nabla \phi\right|^2-\lm<\phi^2>_\de+\lm <\phi>_\de^2}{<\phi^2>_\de},\quad \phi\in H.
$$
 In particular it is not difficult to prove that the first eigenvalue can be characterized as
follows
$$
\tau_1=\inf\{Q(\phi)\,|\, \phi\in H\setminus\{0\}\}.
$$

At this point we argue by contradiction and assume that \eqref{7.1} admits a non trivial solution.
Hence, in particular, $\tau_1\leq 0$ and we readily conclude that
$$
\tau_0:=\inf\{Q_0(\phi)\,|\, \phi\in H\setminus\{0\}\}\leq 0,\quad\mbox{where}\quad
Q_0(\phi)=\frac{\mathcal{L}_0(\phi,\phi)}{<\phi^2>_\de}
$$

and

$$
\mathcal{L}_0(\phi,\psi)= \int\limits_{\om} \left(\nabla \phi\cdot\nabla \psi\right)
-\lm\int\limits_{\om} \de(u)\phi\psi,\;(\phi,\psi)\in H\times H.
$$

Clearly $\tau_0$ is attained by a simple and positive eigenfunction $\vp_0$ which satisfies

\bel\label{7}
\left\{
  \begin{array}{ll}
    \dsp -\Delta\vp_0 -\lm \de(u)\vp_0 =\tau_0 \de(u)\vp_0 & \mbox{in}\ \ {\om} \\
    \vp_0 =0 & \mbox{on}\ \ \pa{\om}.
  \end{array}
\right.
\eel

Let us recall that we have obtained solutions for $P(\lm,{\om})$ as solutions of $Q(\mu,{\om})$
in the form $u=u_{\a,{\mu},c}$, for some $\mu=\mu(\a)\leq \ov{\mu}_\a$ whose value of $\lm=\lm(\mu,\a,c)$
was then estimated as a
function of $\a$. Therefore, at this point, it is more convenient to look at the linearized problem
in the other way, that is, to go back to $\mu=\lm\left(\int_{\om} e^u\right)^{-1}$.
Hence, let us observe that for a generic value $\mu\leq \ov{\mu}_\a$ \rife{7} takes the form
\bel\label{9}
\left\{
  \begin{array}{ll}
    \dsp -\Delta\vp_0 -\mu K_{\a,\mu,c} \vp_0  =\nu_0 K_{\a,\mu,c} \vp_0 & \mbox{in}\ \ {\om} \\
    \vp_0 =0 & \mbox{on}\ \ \pa{\om},
  \end{array}
\right.
\eel
where
$$
K_{\a,\mu,c}=e^{u_{\a,\mu,c}}\quad\mbox{and}\quad \nu_0=\mu\frac{\tau_0}{\lm}\leq 0.
$$

\brm{\it
Of course, the assertion about the positivity of the first eigenvalues corresponds to the positivity of
$\tau_1$ and $\nu_0$ respectively. Therefore that part of the statement will be automatically proved once we get the desired
contradiction.}
\erm

Since also the linearized problem \eqref{7.1} is rotational, translational and dilation invariant, by arguing
exactly as in the proof of Theorem \ref{t1-intro} we can assume without loss of generality that
$$
\Omega_{\rho,c}:=\{\rho^2 x^2+y^2\leq c\}\subset\Omega\subset\{\rho^2 x^2+y^2\leq 1\}=:\Omega_\rho.
$$

We observe that, by defining
$$
K_{\a,\mu,c}^{(-)}:=e^{v_{\a,\ga^{-}_{\a,c}(\mu),c}}=\left\{
                                                       \begin{array}{ll}
                                                         \left(\frac{1+\ga^{-}_{\a,c}(\mu)^2}{1+\frac{\ga^{-}_{\a,c}(\mu)^2}{c}(\a^2x^2+y^2)}\right)^2 & \mbox{$(x,y)\in\Omega_{\a,c}$} \\
                                                         1 & \mbox{$(x,y)\in\Omega\setminus\Omega_{\a,c}$,}
                                                       \end{array}
                                                     \right.
$$
$$
K_{\a,\mu}^{(+)}:=e^{v_{\a,\ga^{+}_\a(\mu)}}=\left(\frac{1+\ga^{+}_\a(\mu)^2}{1+\ga^{+}_\a(\mu)^2(\a^2x^2+y^2)}\right)^2,\qquad (x,y)\in\Omega_{\a}
$$
we have
$$
K_{\a,\mu,c}^{(-)}\leq K_{\a,\mu,c}\leq K_{\a,\mu}^{(+)}\qquad \textnormal{for any }(x,y)\in\Omega.
$$

In particular, since
$$
K_{\a,\mu}^{(+)}\leq (1+\ga_\a^{+}(\mu)^2)^2\quad\mbox{and}\quad
1\leq K_{\a,\mu,c}^{(-)}\leq (1+\ga_{\a,c}^{-}(\mu)^2)^2\quad \mbox{in}\quad \om,
$$
and
\beq\label{140213.0}
\om\subset T_\a:=\{(x,y)\in\rdue\,|\,|\,x|\leq (\a)^{-1},\;|\,y|\leq 1\},
\eeq

then, by using the fact that
$$
\nu_0=\inf\left\{
\left.\frac{\int\limits_{\om} \left|\nabla \vp\right|^2\,dx-\mu\int\limits_{\om} K_{\a,\mu} \vp^2\,dx}
{\int\limits_{\om} K_{\a,\mu}\vp^2\,dx}\;\right|
\,\vp\in H\right\}\leq0,
$$
it is not difficult to check that, for any $\mu\leq \ov{\mu}_\a=\frac{(1+\a^2)^2}{2}$, the following inequality
holds:
\beq\label{13}
\inf\left\{
\left.\frac{\int\limits_{T_\a} \left|\nabla \vp\right|^2\,dx-\mu(1+\ga_\a^{+}(\mu)^2)^2\int\limits_{T_\a} \vp^2\,dx}
{\int\limits_{T_\a} \vp^2\,dx}\;\right| \,\vp\in H\right\}\leq 0.
\eeq

Hence, there exists $\overline{\mu}_0\leq 0$ such that,
putting $\sg=\sg(\mu,\a)=\mu(1+\ga_\a^{+}(\mu)^2)^2+\overline{\mu}_0$, there exists a weak solution $\phi_0\in H$ of
\bel\label{14}
\left\{
  \begin{array}{ll}
    \dsp -\Delta\phi_0 -\sg\phi_0 =0 & \mbox{in}\ \ T_\a,\\
    \dsp \phi_0=0 & \mbox{on}\quad  \pa T_\a.
  \end{array}
\right.
\eel
It is well known that the minimal eigenvalue $\sg_{\it min}$ of \eqref{14} satisfies
$\sg_{\it min}=\frac{\pi^2}{4}\a^2+\frac{\pi^2}{4}>2(1+\a^2)$ and we conclude that
\beq\label{16}
2(1+\a^2)\leq\sg(\mu,\a)=\mu(1+\ga_\a^{+}(\mu)^2)^2+\overline{\mu}_0.
\eeq

Next, since $\underline{\a}_*(c)<\tfrac{1}{2\sqrt{10}}$, it is not difficult to check that
$\sg=\sg(\mu,\a)$ satisfies
$$
\sg(\mu,\a)\leq 1,
$$
for any $\a\leq\underline{\a}_*(c)$, which is of course a contradiction to \eqref{16}. This fact concludes the first
part of the proof. As for the second one it can be derived by arguing as above with some minor changes as in the
proof Theorem \ref{t1-intro}(b).\fineproof

\bigskip

\section{A multiplicity result}\label{sec:mp}

This section is devoted to the proof of Theorem \ref{mp-intro}.

\proof[Proof of Theorem \ref{mp-intro}]

\emph{(a).}\;\; Let us fix $\lm\in(8\pi,\lm_{a,c})\setminus 8\pi\enne$,
then there exists $k\in\enne^*$ such that $\lm\in(8k\pi,8(k+1)\pi)$.
Let us fix now $k$ distinct points, $x_1,\ldots,x_k$,
in the interior of $\Omega_{\a,\beta_-}=\{\a^2x^2+y^2\leq \beta_-^2\}$.
Next we fix $\bar d>0$ such that $\dist(x_i,x_j)>4\bar d$ for any $i\neq j$
and such that $\dist(x_i,\pa\om_{\a,\beta_-})>2\bar d$ for any $i\in\{1,\ldots,k\}$.

Following \cite{dj} we introduce some notations. For $d\in(0,\bar d)$ we consider
a smooth non-decreasing cut-off function $\chi_d:[0,+\infty)\to\R$ satisfying the following properties:

$$\left\{
    \begin{array}{ll}
      \chi_d(t)=t & \hbox{for $t\in[0,d]$} \\
      \chi_d(t)=2d & \hbox{for $t\geq 2d$} \\
      \chi_d(t)\in[d,2d] & \hbox{for $t\in[d,2d]$.}
    \end{array}
  \right.
$$
Then, given $\mu>0$, we define the function $\varphi_{\mu,d}\in H^1_0(\Omega)$ by
$$
\varphi_{\mu,d}(y)=\left\{
                          \begin{array}{ll}
                            \log\,\sum_{j=1}^k \frac 1k
\left(\frac{8\mu^2}{(1+\mu^2\chi^2_d(|y,x_j|))^2}\right)-\log\left(\frac{8\mu^2}{(1+4d^2\mu^2)^2}\right)
& \hbox{$y\in\Omega_{\a,\beta_-}$} \\
                            0 & \hbox{$y\in\Omega\setminus\Omega_{\a,\beta_-}$.}
                          \end{array}
                        \right.
$$

By arguing exactly as in Section 5 of \cite{dj} we have
$$
F_{\lm}(\varphi_{\mu,d})\leq(16k\pi-2\lm+o_d(1))\ln(\mu)+O(1)+C_d
$$
where $C_d$ is a constant independent of $\mu$ and $o_d(1)\to 0$ as $d\to 0$.\\
Then, there exist $d_0$ sufficiently small and $\mu_0$ sufficiently large such that
$$
F_{\lm}(\varphi_{\mu_0,d_0})<F_{\lm}(\ul)-1.
$$
Next we define
$$
\mathcal{D}=\{\gamma:[0,1]\to H^1_0(\Omega)\,:\,\gamma\textnormal{ is continuous, $\gamma(0)=\ul$, $\gamma(1)=\varphi_{\mu_0,d_0}$}\}
$$
and, for any $\eta\in(8k\pi,8(k+1)\pi)\cap(8\pi,\lm_{a,c})$, we set
$$
c_\eta=\inf\limits_{\gamma\in\mathcal{D}}\max\limits_{s\in[0,1]}F_\eta(\gamma(s)).
$$
Since $\ul$ is a strict local minimum for $F_\lm$, there exists $\eps_\lm>0$ such that $c_\lm\geq F_\lm(\ul)+\eps_\lm$.
Besides, since $F_\lm$ is continuous and the branch $\mathcal{G}_{\a,c}$ is smooth,
we have that a bound on the min-max levels applies uniformly in a small neighborhood of $\lm$.
More precisely the following straightforward fact holds true.
\ble\label{lemma:lambda0}
There exists $\lm_0>0$ sufficiently small such that
$$
[\lm-\lm_0,\lm+\lm_0]\subset(8k\pi,8(k+1)\pi)\cap(8\pi,\lm_{\a,c})
$$
and for any $\eta\in[\lm-\lm_0,\lm+\lm_0]$ we have $F_\eta(\varphi_{\mu_0,d_0})\leq F_\eta(\ul)-\frac12$ and
$$
c_\eta\geq F_\lm(\ul)+\frac34\eps_\lm\geq F_\eta(\ul)+\frac12\eps_\lm.
$$
\ele

If $\eta$, $\eta'\in(\lm-\lm_0,\lm+\lm_0)$, $\eta\leq \eta'$, then
$
\frac{F_\eta}{\eta}-\frac{F_\eta'}{\eta'}=\frac12(\frac{1}{\eta}-\frac{1}{\eta'})\ino|\grad u|^2\geq 0,
$
whence
\bel\label{nonincreasing}
\frac{c_\eta}{\eta}\geq\frac{c_{\eta'}}{\eta'}.
\eel
Therefore we have that the function $\eta\mapsto \frac{c_\eta}{\eta}$ is non-increasing and in turn
differentiable a.e. in $(\lm-\lm_0,\lm+\lm_0)$.
Set
$$
\Lambda=\{\eta\in(\lm-\lm_0,\lm+\lm_0)\,|\,\frac{c_\eta}{\eta}\textrm{ is differentiable at $\eta$}\}.
$$

\ble\label{lemma:denso}
$c_\eta$ is achieved by a critical point $\vl$ of $F_\eta$ provided that $\eta\in\Lambda$.
\ele
\proof The proof is a step by step adaptation of the arguments of Lemma 3.2 of
\cite{DJLW}  where, with respect to their notations, we have just to choose $\dt<\frac14\eps_{\lm}$.
\fineproof

\bigskip

Finally we state a (well known) compactness result for sequence of solutions of
$P(\lm_n,\Omega)$.
\ble\label{lemma:compattezza}
Let $\lm_n\to\lm$ and let $v^{\scp(\lm_n)}\in H^1_0(\Omega)$ be a solution of $P(\lm_n,\Omega)$.
If $\lm\notin 8\pi\enne$, then $v^{\scp(\lm_n)}$ admits a subsequence which converges smoothly
to a solution $\vl$ of $P(\lm,\Omega)$.
\ele
\proof In view of Lemma 2.1 in \cite{CCL} $v^{\scp(\lm_n)}$ is uniformly bounded
in a \un{fixed} neighborhood of the boundary.
Hence the conclusion is a straightforward and well known consequence of the Brezis-Merle \cite{bm}
concentration-compactness result as completed by Li and Shafrir \cite{ls}.
\fineproof

\bigskip

Now we are able to conclude the proof of Theorem \ref{mp-intro}(a).
Indeed the thesis is an easy consequence of Lemmas \ref{lemma:denso} and \ref{lemma:compattezza},
noticing that the solution $\vl$, obtained by this procedure,
does not coincide with $\ul$, because by Lemma \ref{lemma:lambda0} $F_\lm(\vl)> F_\lm(\ul)$.

\

\emph{(b).}\;\;This part can be proved exactly as the previous one.
\fineproof

\section{A refined estimate for solutions on $\mathcal{G}_{\a,1}$}\label{s5}
Let $\mathcal{G}_{\a,c},\mathcal{G}_\textnormal{\tiny{$N$}}$ denote the branches of parameter-solutions pairs
of $P(\lm,\om)$ found in Theorem \ref{t1-intro}.
As a consequence of Theorem \ref{unique:1-intro} and Proposition
\ref{pr2}  we obtain the following:

\bpr\label{pr3:I}
Let $\ov{\lm}\geq 8\pi$, $\widetilde{\a}_{1}$ and $\widetilde{N}$ be as in Theorem \ref{unique:1-intro}. Let
either
$\mathcal{G}^{(\ov{\lm})}=\{(\lm,\ul)\in\mathcal{G}_{\a,c}\,:\,\lm\in[0,\ov{\lm})\}$ or
$\mathcal{G}^{(\ov{\lm})}=\{(\lm,\ul)\in\mathcal{G}_\textnormal{\tiny{$N$}}\,:\,\lm\in[0,\ov{\lm})\}$
denote that part of $\mathcal{G}_{\a,c},\mathcal{G}_\textnormal{\tiny{$N$}}$
 with $\lm\in[0,\ov{\lm})$, $\a\in(0,\widetilde{\a}_{1}]$
and $N\geq \widetilde{N}$ respectively. Then the energy function
\beq\label{endef:II}
\widehat{E}(\lm):= \mathcal{E}(\de(\ul)),\quad \ul \in \mathcal{G}^{(\ov{\lm})},
\eeq
is a monodrome and smooth function of $\lm\in[0,\ov{\lm})$.
\epr
\proof
By using the explicit bounds \rife{susu} and the fact that
$$
\mathcal{E}(\de(\ul))=\frac{1}{2\lm}\ino\de(\ul)\ul,
$$
then it is straightforward to show that the energy of any solution lying on $\mathcal{G}^{(\ov{\lm})}$ is
uniformly bounded
from above by a suitable value $\ov{E}$, which we can assume without loss of generality to be larger than $1$.
Therefore Theorem \ref{unique:1-intro} applies and we see that
$\mathcal{E}(\de(\ul))$ is monodrome as a function of $\lm\in[0,\ov{\lm})$
and consequently $\widehat{E}(\lm)$ is well defined. At this point Proposition \ref{pr2} implies that
it is smooth as well, see also Remark \ref{newrmsmooth}.
\fineproof

\bigskip

Our next aim is to improve Proposition \ref{pr3:I} in case $\om=\om_\a$ to come up with a unique solution of
$P(\lm,\om_\a)$ at fixed energy. Indeed, this is the content of Theorem \ref{unique:2-intro} whose proof is
the main aim of this section. To achieve this goal we have to pay a price in
terms of a smallness assumption on the energy and indeed we will obtain this result
by using Theorem \ref{unique:1-intro} and the expansion of solutions as functions of $\a$.
Actually, we first need a more precise formula about the explicit form of solutions
of $P(\lm,\om_\a)$ lying on $\mathcal{G}_{\a,1}$, as claimed in \rife{sol:II-intro} of
Theorem \ref{thm:261112-intro}.
By using these expansions we will be able to calculate explicitly, at least for small $\a$,
their energy as a function of $\lm$ and then prove that $\widehat{E}$ is monotone. It turns out
that this is enough to prove uniqueness of solutions with fixed energy.
Actually we also provide another proof (still by using the sub-supersolutions method)
of the existence of solutions for $P(\lm,\om_\a)$.

\bigskip

{\bf The Proof of Theorem \ref{thm:261112-intro}.}\\
As above, the notation $\mbox{O}(\a^m)$, $m\in\enne$ will be used in the rest of this proof to denote various quantities
uniformly bounded by $C_m\a^m$ with $C_m>0$ a suitable constant depending only on $\ov{\lm}$.\bigskip

Let us first seek solutions
$v_\a$ of $Q(\mu_0\rho,\om_\a)$ in the form

\beq\label{021212.1}
v_\a=\a \phi_0+\a^2\phi_{0,1},\quad \phi_0,\phi_{0,1}\in C^2(\om_\a)\cap C^0(\ov{\om_\a}),
\eeq
with the additional constraints

$$
0\leq \|\phi_0\|_\i\leq M_0,\quad 0\leq \|\phi_{0,1}\|_\i\leq M_1.
$$

Since $v_\a$ must satisfy $-\Delta v= \mu_0 \a{\dsp e^v}$ then $\phi_0$ and $\phi_{0,1}$ should be solutions of

\beq\label{191112.1}
\graf{
-\Delta \phi_0  =\mu_0 & \mbox{in}\quad \om_\a\\
\phi_0 =0 & \mbox{on}\quad \pa\om_\a
}
\eeq
and
\beq\label{231012.1.II}
\graf{
-\Delta \phi_{0,1}=\mu_0\a^{-1}\left(e^{\a\phi_0}e^{\a^2\phi_{0,1}}-1\right) & \mbox{in}\quad \om_\a\\
\phi_{0,1}=0 & \mbox{on}\quad \pa\om_\a
}
\eeq
respectively. Therefore the explicit expression of $\phi_0$ is easily derived to be
\beq\label{021212.2}
\phi_0(x,y;\a)=\frac{\mu_0}{2(1+\a^2)}\left(1-(\a^2x^2+y^2)\right),\quad (x,y)\in \om_\a.
\eeq
Please observe that the function $\phi_0(x,y;\lm,\a)$ as defined in \rife{phi0-intro}
will be recognized to be $\phi_0(x,y;\a)$ where $\mu_0=\mu_0(\lm,\a)$.\\

Clearly
$$
\|\phi_0\|_\i=\frac{\mu_0}{2(1+\a^2)},
$$
and therefore,  in particular we have
\beq\label{231012.2.II}
\fo t_0>1\, \E\,\a_1=\a_1(t_0)>0\,:\,
e^{\a\phi_0}\leq e^{\frac{\mu_0\a}{2(1+\a^2)}}<1+t_0\frac{\mu_0\a}{2},\quad \fo\a<\rho_1,
\eeq
the last inequality being a trivial consequence of the convexity of $e^{\frac{\mu_0 s}{2(1+s^2)}}$ in a right neighborhood of $s=0$.\\
Our next aim is to use the sub-supersolutions method to obtain solutions for \rife{231012.1.II}. Let us define
$$
f(t;\phi_0):=e^{\a\phi_0}e^{\a^2t},\quad t\geq 0,
$$

so that, in particular, we have
\beq\label{231012.3.II}
\fo t_1>1\,\E\,\a_2>0\,:\,e^{\a^2t}<1+t_1 \a^2 t,\quad \fo\a<\a_{2},
\eeq
with $\a_2$ depending on  $t_1$.
By using \rife{231012.2.II} and \rife{231012.3.II} we conclude that
$$
f(t;\phi_0)\leq \left(1+t_0\frac{\mu_0\a}{2}\right)\left(1+t_1 \a^2 t\right),\quad\fo\a<\min\{\a_1,\a_2\}.
$$

Hence, by setting
$$
A_{+}=1+t_0\frac{\mu_0\a}{2},
$$
we see that a supersolution $\phi_{+}$ for \rife{231012.1.II} will be obtained whenever we will be able
to solve the differential problem

\beq\label{231012.4.II}
\graf{
-\Delta \phi_{+}\geq t_0\frac{\mu_0^2}{2} + t_1\mu_0 A_{+} \a \phi_{+} & \mbox{in}\quad \om_\a\\
\phi_{+}\geq 0 & \mbox{on}\quad \pa\om_\a\\
0\leq\phi_{+}\leq M_1 & \mbox{in}\quad \om_\a.
}
\eeq

Let us define
$$
\phi_+(x,y)=\frac{C_{+}}{2(1+\a^2)}\left(1-(\a^2x^2+y^2)\right),\quad (x,y)\in \om_\a,
$$
with $C_+>0$, so that the differential inequality in \rife{231012.4.II} yields
$$
-\Delta \phi_+=C_+=\frac{C_+}{2}+\frac{C_+}{2}=\frac{C_+}{2}+(1+\a^2)\|\phi_+\|_\i\geq t_0\frac{\mu_0^2}{2} + t_1\mu_0 A_{+} \a \phi_{+}.
$$

Therefore \rife{231012.4.II} will be satisfied whenever we can choose $C_+$ such that the following inequalities are verified
\beq\label{231012.5.II}
\graf{
C_+ &\geq t_0\mu_0^2\\
(1+\a^2)&\geq t_1\mu_0 A_{+} \a\\
C_+&\leq 2(1+\a^2)M_1.
}
\eeq

We first impose
$$
C_+=2M_1,
$$
so that the third inequality in \rife{231012.5.II} is automatically satisfied and then substitute it in the first inequality,
to obtain
\beq\label{021212.3}
\mu_0^2\leq \min\left\{\frac{2M_1}{t_0},\left(4M_0\right)^2\right\}={\frac{2M_1}{t_0}},
\;\mbox{for any}\;M_0\;\mbox{large enough}.
\eeq

We conclude in particular that the second inequality is trivially satisfied for any $\a$ small enough. At this point Theorem \ref{t0}
shows that there exists a solution $v_\a$ of $Q(\mu_0\a,\om_\a)$ taking the form \rife{021212.1},
where $\phi_0$ is defined as in \rife{021212.2}
and $0\leq \phi_{0,1}\leq M_1$ with the constraint \rife{021212.3}.\\

Our next aim is to show that $\fo\ov{\lm}\geq 8\pi$ we can find $\a_0$ small enough such that $\fo\a<\a_0$ and for any $\lm<\ov{\lm}$ we can
choose $\mu_0$ in such a way that $v_\a$ is a solution of $P(\lm,\om_\a)$. Indeed, we have
\beq\label{vip}
\lm=\lm_0(\mu_0,\a):=\mu_0\a\int\limits_{\om_\a} e^{v_\a}=\pi\mu_0+f_0(\mu_0,\a),\;\;\mbox{where}\;\;
|f_0(\mu_0,\a)|\leq C_{\scp M_1}\a^2,
\eeq

where $\lm$ is a fixed value in the range of $\lm_0$ and we have used $\|\phi_{0,1}\|_\i\leq M_1$ and
$$
\int\limits_{\om_\a} e^{\a\phi_0}=\left(1+\a^2\right)\frac{2\pi}{\mu_0\a^3}\left(e^{\frac{\mu_0\a^2}{2(1+\a^2)}}-1\right).
$$

\ble\label{intermedio}
$\lm_0(\mu_0,\a)$ is smooth.
\ele
\proof
It is straightforward to check that the energy of these solutions $v_\a$ is uniformly bounded
from above by a suitable positive number $\ov{E}$ (possibly depending on $M_1$ and $\ov{\lm}$)
which we can assume without loss of generality to be larger than 1.
Therefore Theorem \ref{unique:1-intro} shows that they must coincide
with some subset of the branch $\mathcal{G}^{(\ov{\lm})}$ (see Proposition \ref{pr3:I}).
We can use Proposition \ref{pr2} at this point and conclude that $\lm_0(\mu_0,\a)$ is smooth as a function of
$\mu_0$. At this point the (joint) regularity of $\lm_0(\mu_0,\a)$ as a function of $\mu_0$ and $\a$
is derived by a conformal transplantation on the unit disk, classical representation formulas
for derivatives of  Riemann maps (see for example \cite{pom}) and standard elliptic theory.
\fineproof

\bigskip

Hence, in particular
we can always choose $\mu_0$ and $\a_0$ such that $\fo\a<\a_0$ we have (see \rife{021212.3})
$$
[0,\ov{\lm})\subset\lm_0\left(\left[0,2\sqrt{\frac{M_1}{t_0}}\right),\a\right),
$$
and since $\lm_0(\mu_0,\a)$ is also continuous, we finally obtain the desired solution for any $\lm<\ov{\lm}$.

At this point, let us fix a positive value $\lm<\ov{\lm}$ for which we seek an approximate solution
$u_{\scp \lm}$ of $P(\lm,\om_\a)$. As a consequence of \rife{vip} we have

\beq\label{mufix}
\mu_0=\mu_0(\lm,\a)=\frac{\lm}{\pi}+\mbox{O}(\a^2),
\eeq

and then

\beq\label{ulfix}
u_{\scp \lm}:=\a \phi_0+\a^2\phi_{0,1}=\frac{\a\mu_0}{2\pi(1+\a^2)}\left(1-(\a^2x^2+y^2)\right)\left(1+\mbox{O}(\a)\right)=
\frac{\a\lm}{2\pi}\left(1-(\a^2x^2+y^2)\right)\left(1+\mbox{O}(\a)\right),
\eeq

is a solution for  $P(\lm,\om_\a)$, as desired.\\

\brm
However, by using \rife{lm0-intro}, \rife{mufix} and \rife{ulfix}, a straightforward explicit evaluation shows that
$$
\mathcal{E}(\de_{\scp \lm})=\frac{1}{2}\int\limits_{\om_\a} \de_{\scp \lm} G_\a[\de_{\scp \lm}]=
\frac{1}{2\lm}\int\limits_{\om_\a} \de_{\scp \lm} \ull=\frac{\mu_0\a}{2\lm^2}\int\limits_{\om_\a}\e{\ull} \ull=
$$
$$
\frac{\lm\a+\mbox{O}(\a^3)}{2\pi\lm^2}\int\limits_{\om_\a}\e{\ull} \ull=\frac{\a}{8\pi}\left(1+\mbox{O}(\a)\right),
$$
see Remark \ref{rem6.1}. Therefore, as far as we are interested in the monotonicity of
$\widehat{E}(\lm)$, we see that the first order expansion is not enough to our purpose.
\erm

Hence we make a further step to come up with an expansion of $\mathcal{E}$ at order $\a^2$.
Let $\phi_{0,1}$ be the solution
of \rife{231012.1.II} determined above, we write it as
$$
\phi_{0,1}=\phi_1+\a\phi_2,
$$
so that, if $\phi_1$ is the unique solution of
\beq\label{011112.0}
\graf{
-\Delta \phi_1 =\mu_0\phi_0=\mu_0^2\psi_0 & \mbox{in}\quad \om_\a\\
\phi_1 =0 & \mbox{on}\quad \pa\om_\a
}
\eeq
(see \rife{phi0-intro}-\rife{080213.1-intro})
then by definition $\phi_2$ is a solution for
\beq\label{011112.1}
\graf{
-\Delta \phi_2=\mu_0\a^{-1}\left(e^{\a\phi_0}e^{\a^2\phi_{0,1}}-1-\phi_0\right) & \mbox{in}\quad \om_\a\\
\phi_2=0 & \mbox{on}\quad \pa\om_\a
}
\eeq
and it is not difficult to check that it also satisfies $\|\phi_2\|\leq M_2$, for a suitable
$M_2$ depending only  $M_0$ and $M_1$.\\
At this point standard elliptic estimates to be used together with the maximum principle show
that $\{\phi_0,\phi_1,\phi_2\}\subset C^{2}_0(\om)$ and
$\|D^{(k)}_\lm\phi_0\|_{\scp C^{2}_0(\om)}+\|D^{(k)}_\lm\phi_1\|_{\scp C^{2}_0(\om)}+
\|D^{(k)}_\lm\phi_2\|_{\scp C^{2}_0(\om)}\leq \ov{M}_k$
for suitable constants $\ov{M}_k>0$ depending only on $M_0$, $M_1$, $M_2$, that is, depending only on $\ov{\lm}$.\\

Let $\lm_0=\lm_0(\mu_0,\a)$ as defined in \rife{vip} above. In view of Lemma \ref{intermedio},
we can expand $\lm_0$ at second order in $\a$,
$$
\lm_0(\mu_0,\a):=\mu_0\a\int\limits_{\om_\a} e^{v_\a}=
\mu_0\a\int\limits_{\om_\a} (1+ \a\phi_0
+\mbox{O}(\a^2))=
$$
$$
\pi\mu_0+\frac{\pi\mu_0^2\a}{4(1+\a^2)}+\mbox{O}(\a^2)=\pi\mu_0+\frac{\pi\mu_0^2\a}{4}+\mbox{O}(\a^2).
$$

Hence, for a fixed value $\lm$ in the range of $\lm_0$ we can use the implicit function theorem to obtain the
inverse expansion up to order $\a^2$, that is
$$
\lm=\pi\mu_0+\frac{\pi\mu_0^2\a}{4}+\mbox{O}(\a^2),\quad\mu_0=\frac{\lm}{\pi}-\frac{\lm^2}{4\pi^2}\a+\mbox{O}(\a^2),
$$
and \rife{mu0-intro}-\rife{mu0-intro1} follows immediately.\\
This observation concludes the proof.\fineproof

\bigskip

{\bf The Proof of Theorem \ref{unique:2-intro}}\\
The notation $\mbox{O}(\a^m)$, $m\in\enne$ will be used in the rest of this proof to denote various quantities
uniformly bounded by $C_m\a^m$ with $C_m>0$ a suitable constant possibly depending on $\ov{\lm}$ and on
the constants $\ov{M}_k$, $k=1,2,3$ as obtained in Theorem \ref{thm:261112-intro}.\\

By using \rife{lm0-intro} above and Theorem \ref{thm:261112-intro} we obtain the Taylor expansion

$$
\mathcal{E}(\de_{\scp \lm})=\frac{1}{2}\int\limits_{\om_\a} \de_{\scp \lm} G_\a[\de_{\scp \lm}]=
\frac{1}{2\lm}\int\limits_{\om_\a} \de_{\scp \lm} \ull=\frac{\mu_0\a}{2\lm^2}
\int\limits_{\om_\a}\e{\ull} \ull=
$$

$$
\frac{\mu_0\a}{2\lm^2}\int\limits_{\om_\a}\e{\ull} \ull=
\frac{\mu_0\a}{2\lm^2}\int\limits_{\om_\a}
(1+ \a \phi_0 +   \mbox{O}(\a^2))(\a \phi_0 +  \a^2\phi_1 +\mbox{O}(\a^3))=
$$

$$
\frac{\mu_0\a}{2\lm^2}\int\limits_{\om_\a}(\a \phi_0 + \a^2 \phi^2_0 + \a^2\phi_1 +\mbox{O}(\a^3))
=\frac{\mu_0\a}{2\lm^2}\left[\frac{\pi\mu_0}{4(1+\a^2)}+\frac{\pi\mu^2_0\a}{12(1+\a^2)^2}
+\frac{\pi\mu^2_0\a}{12(1+\a^2)^2}+\mbox{O}(\a^2) \right],
$$
where we have used the fact that
\beq\label{011112.3}
\int\limits_{\om_\a}\a^2\phi_1=\frac{\pi\mu^2_0\a}{12(1+\a^2)^2},
\eeq
which can be obtained by using the explicit expression of $\phi_0$ in \rife{phi0-intro}
together with the fact that $\phi_1$
solves \rife{011112.0}, see the Appendix \ref{subs:1} below for further details.\\

Hence, by using Proposition \ref{pr3:I} and \rife{mu0-intro}-\rife{mu0-intro1} and \rife{regular-intro}, we have

$$
\widehat{E}(\lm):=\mathcal{E}(\de_{\scp \lm})=\frac{\pi \mu_0^2\a}{8\lm^2}+\frac{\pi \mu_0^3\a^2}{12\lm^2}+\mbox{O}(\a^3)=
\frac{\pi \a}{8\lm^2}\left(\frac{\lm^2}{\pi^2}-\frac{\lm^3}{2\pi^3}\a+\mbox{O}(\a^2)\right)+
$$
$$
\frac{\pi \a^2}{12\lm^2}\frac{\lm^3}{\pi^3}+\mbox{O}(\a^3)=
\frac{\a}{8\pi}+\frac{\a^2}{48\pi^2}\lm+\mbox{O}(\a^3).
$$

In particular we conclude that
\beq\label{231112.3.8}
\widehat{E}(\lm)=\frac{\a}{8\pi}+\frac{\a^2}{48\pi^2}\lm+\mbox{O}(\a^3),
\eeq

and, in view of \rife{mu0-intro}-\rife{mu0-intro1} and \rife{regular-intro},
\beq\label{231112.3.a}
\frac{d}{d\lm}\widehat{E}(\lm)=\frac{\a^2}{48\pi^2}+\mbox{O}(\a^3),
\eeq
$$
\frac{d^2}{d\lm^2}\widehat{E}(\lm)=\mbox{O}(\a^3).
$$

At this point \rife{231112.3.8} shows that we may
restrict the domain of $\widehat{E}$ to the preimages of
$E\in\left[\frac{\a}{8\pi},\widehat{E}_\a\right]$. Then \rife{231112.3.a} implies that $\widehat{E}(\lm)$
is monotonic increasing there. Hence the preimage of
$\left[\frac{\a}{8\pi},\widehat{E}_\a\right]$ is exactly $[0,\widehat{\lm}_\a]$
and the uniqueness of $u_{\scp \lm}$ as a function of $\lm$ implies that the equation
$\mathcal{E}(\de(u_{\scp \widehat{\lm}(E)}))=E$ defines $\widehat{\lm}(E)$
as a monotonic increasing function of $E$ in $\left[\frac{\a}{8\pi},\widehat{E}_\a\right]$.
Therefore, we can use \rife{231112.3.8} and \rife{231112.3.a} together with the implicit function theorem to
take the inverse up to order $\a^2$, that is
$$
\widehat{\lm}(E)=\frac{48\pi^2}{\a^2}\left(E-\frac{\a}{8\pi}\right)+\mbox{O}(\a),
$$
and then conclude that
$$
\frac{d}{d E}\widehat{\lm}(E)=\frac{48\pi^2}{\a^2}+\mbox{O}(\a),
$$
and
$$
\frac{d^2}{d E^2}\widehat{\lm}(E)=\mbox{O}(\a).
$$
\fineproof

\section{The Entropy is concave in $E\in\left[\frac{\a}{8\pi},\widehat{E}_\a\right]$.}\label{sec:entropy}

Let us recall that according to definition \ref{dedef} the density corresponding to a solution
$\ull$ of $P(\lm,\om_\a)$ is defined to be
$$
\de_{\scp \lm}\equiv\de(\ull):=\dfrac{\e{\ull}}{\int\limits_{\om_\a} \e{\ull}}.
$$
As usual $\mathcal{G}_{\a,1}$ denotes the branch of solutions obtained in Theorem \ref{t1-intro}(a).\\

When evaluated on $(\lm,\ull)\in\mathcal{G}_{\a,1}$, of course $\mathcal{S}(\de(\ull))$ yields a function
of $\lm$ defined in principle on $\lm\in[0,\lm_{\a,1}]$.
Then we can use MVP-(iv), that is, the fact that any entropy maximizer (at fixed $E$) of the MVP satisfies
$P(\lm,\om)$ (for a certain unknown value $\lm$). But then we can observe that
Theorem \ref{unique:2-intro} states that there exists one and only one solution of
$P(\lm,\om)$ with $\lm=\widehat{\lm}(E)$ such that
the energy is exactly $E$, $\widehat{E}(\lm)=E$, whenever $E\in\left[\frac{\a}{8\pi},\widehat{E}_\a\right]$.\\
Therefore we conclude that indeed $S(E)\equiv \mathcal{S}(\de(\ull))\left.\right|_{\lm=\widehat{\lm}(E)}$ in
$\left[\frac{\a}{8\pi},\widehat{E}_\a\right]$.
Hence, when evaluated on those densities $\de_{\scp \widehat{\lm}(E)}$ as obtained in Theorem \ref{unique:2-intro},
we have
$$
S(E)\equiv\mathcal{S}(\de_{\scp \widehat{\lm}(E)})=-2E \widehat{\lm}(E)+
\log\left(\,\int\limits_{\om_\a} e^{u_{\scp \widehat{\lm}(E)}}\right),\quad E\in\left[\frac{\a}{8\pi},\widehat{E}_\a\right].
$$

In particular, in view of Theorem \ref{unique:2-intro} we can set

$$
\dot{u}=\frac{du_{\scp \widehat{\lm}(E)}}{d E},\quad\mbox{and}\quad
\ddot{u}=\frac{d^2u_{\scp \widehat{\lm}(E)}}{d E^2},
$$
to obtain
$$
\frac{d S(E)}{d E}=-2\widehat{\lm}(E)-2E \frac{d \widehat{\lm}(E)}{d E}+
\int\limits_{\om_\a}\de_{\scp \widehat{\lm}(E)}\dot{u},
$$
and then
\beq\label{entropy}
\frac{d^2 S(E)}{d E^2}=-4\frac{d \widehat{\lm}(E)}{d E}-2E \frac{d^2 \widehat{\lm}(E)}{d E^2}+
\int\limits_{\om_\a}\de_{\scp \widehat{\lm}(E)}(\dot{u})^2-
\left(\int\limits_{\om_\a}\de_{\scp \widehat{\lm}(E)}\dot{u}\right)^2+
\int\limits_{\om_\a}\de_{\scp \widehat{\lm}(E)}\ddot{u}.
\eeq

We wish to evaluate $\frac{d^2 S(E)}{d E^2}$ in case $\om=\om_\a$ and
$E\in\left[\frac{\a}{8\pi},\widehat{E}_\a\right]$. Indeed, this is the content of Proposition
\ref{pr:entropy-intro}.

{\bf The Proof of Proposition \ref{pr:entropy-intro}}\\
We are going to evaluate \rife{entropy} by using \rife{mu0-intro}-\rife{mu0-intro1}, Theorem \ref{unique:2-intro} and the
estimates \rife{regular-intro} in Theorem \ref{thm:261112-intro}. Let us set
$$
\dot{\widehat{\lm}}=\frac{d }{d E}\widehat{\lm}(E),\quad \ddot{\widehat{\lm}}=\frac{d^2 }{d E^2}\widehat{\lm}(E),
$$
and
$$
\phi_j^{'}=\frac{d }{d \lm}\,\phi_j,\quad \phi_j^{''}=\frac{d^2 }{d \lm^2}\,\phi_j,\qquad j=0,1,2,
$$

so that, in view of \rife{regular-intro} and
\rife{231112.3-intro}, \rife{231112.3.a-intro}, \rife{231112.3.b-intro} we have
$$
\dot{u}=\frac{d u}{d\lm}\dot{\widehat{\lm}}=
\dot{\widehat{\lm}}\left( \a \phi_0^{'}+\a^2 \phi_1^{'}+\mbox{\rm O}(\a^3)\right),
$$
and
\beq\label{der1}
\ddot{u}=\frac{d^2 u}{d\lm^2}\dot{\widehat{\lm}}^2+\frac{d u}{d\lm}\ddot{\widehat{\lm}}=
\dot{\widehat{\lm}}^2\left( \a \phi_0^{''}+\a^2 \phi_1^{''}+\mbox{\rm O}(\a^3)\right)+
\ddot{\widehat{\lm}}\left( \a \phi_0^{'}+\a^2 \phi_1^{'}+\mbox{\rm O}(\a^3)\right),
\eeq
where the derivatives with respect to $\lm$ will be estimated by using \rife{mu0-intro}-\rife{mu0-intro1}.\\
Hence we can introduce
$$
\ddot{S}_0(E):=\int\limits_{\om_\a}\de(u_{\widehat{\lm}(E)})(\ddot{u}+\dot{u}^2)
-\left(\int\limits_{\om_\a}\de(u_{\widehat{\lm}(E)})\dot{u}\right)^2,
$$
to obtain, after a lengthy evaluation where we use \rife{mu0-intro}-\rife{mu0-intro1} and \rife{der1},
$$
\ddot{S}_0(E)=\frac{(48\pi)^2}{\a^2}
\left[
-\int\limits_{\om_\a}\de(u_{\widehat{\lm}(E)})\psi_0-
\left(\int\limits_{\om_\a}\de(u_{\widehat{\lm}(E)})\psi_0\right)^2+
\pi^2\int\limits_{\om_\a}\de(u_{\widehat{\lm}(E)})\phi_1^{''}
\right]+\mbox{\rm O}\left(\frac{1}{\a}\right).
$$

At this point we can use

\beq\label{asym:1.0}
\int\limits_{\om_\a}\de(u_{\widehat{\lm}(E)})\psi_0=\frac14+\mbox{\rm O}(\a),
\eeq
and
\beq\label{asym:2.0}
\int\limits_{\om_\a}\de(u_{\widehat{\lm}(E)})\phi_1^{''}=\frac{1}{6\pi^2}+\mbox{\rm O}(\a),
\eeq
whose proof is left to Appendix \ref{subs:2}, and \rife{231112.3.a-intro}, \rife{231112.3.b-intro} to obtain
$$
\frac{d^2 S(E)}{d E^2}=-4\dot{\widehat{\lm}}-2E\ddot{\widehat{\lm}}+\ddot{S}_0(E)=
-4\frac{48 \pi^2}{\a^2}+\frac{(48\pi)^2}{\a^2}\left(-\frac14-\frac{1}{16}+\frac{1}{6}\right),
$$
and the conclusion readily follows.
\fineproof

\section{Appendix}

\subsection{ The proof of \rife{011112.3}}\label{subs:1}$\left.\right.$\\
To obtain \rife{011112.3} we multiply $-\Delta \phi_{1}$ by $y^2$ and integrate by parts twice to obtain
$$
-\int\limits_{\om_\a} y^2 \Delta \phi_{1} = - \int\limits_{\pa \om_\a} y^2 \pa_\nu \phi_{1}-
2\int\limits_{\om_\a} \phi_{1}.
$$

Similarly we have

$$
-\int\limits_{\om_\a} \a^2x^2 \Delta \phi_{1} = - \int\limits_{\pa \om_\a}\a^2 x^2 \pa_\nu \phi_{1}-
2\a^2\int\limits_{\om_\a} \phi_{1},
$$
so that we can sum up to obtain
$$
2(1+\a^2)\int\limits_{\om_\a} \phi_{1}= \int\limits_{\om_\a} (\a^2x^2+y^2) \Delta \phi_{1}-
\int\limits_{\pa \om_\a}\pa_\nu \phi_{1}.
$$

Therefore, by using the equation in \rife{011112.0} and the divergence theorem we have

$$
2(1+\a^2)\int\limits_{\om_\a} \phi_{1} = \int\limits_{\om_\a} (-(\a^2x^2+y^2) +1)\mu_0\phi_0,
$$
that is

\beq\label{19.0.1}
\int\limits_{\om_\a} \phi_{1} =(\mu_0)^2\int\limits_{\om_\a}\psi_0^2,
\eeq
and the conclusion follows by a straightforward evaluation based on the explicit expression of $\psi_0$
(see \rife{080213.1-intro}).\fineproof

\subsection{ The proofs of \rife{asym:1.0} and \rife{asym:2.0}}\label{subs:2}$\left.\right.$\\
Concerning \rife{asym:1.0} we just observe that

$$
\int\limits_{\om_\a}\de(u_{\widehat{\lm}(E)})\psi_0=
\int\limits_{\om_\a}\frac{1+\mbox{\rm O}(\a)}{\int\limits_{\om_\a}(1+\mbox{\rm O}(\a))}\psi_0=
\frac{\a}{\pi}(1+\mbox{\rm O}(\a))\int\limits_{\om_\a}\psi_0=
\frac14+\mbox{\rm O}(\a),
$$
where the last equality is obtained by a straightforward evaluation based on the explicit expression of $\psi_0$
(see \rife{080213.1-intro}).\\

Concerning \rife{asym:2.0} we observe as above that
\beq\label{19.0.2}
\int\limits_{\om_\a}\de(u_{\widehat{\lm}(E)})\phi_1^{''}=
\frac{\a}{\pi}(1+\mbox{\rm O}(\a))\int\limits_{\om_\a}\phi_1^{''},
\eeq
and that in view of \rife{011112.0} and \rife{phi0-intro}, then $\phi_1^{''}$ satisfies
\beq\label{19.0.0}
\graf{
-\Delta \phi_1^{''} =(\mu_0\phi_0)^{''}\equiv(\mu_0^2)^{''}\psi_0 & \mbox{in}\quad \om_\a\\
\phi_1^{''} =0 & \mbox{on}\quad \pa\om_\a
}
\eeq

where $\mu_0=\mu_0(\lm,\a)$ (see \rife{mu0-intro}-\rife{mu0-intro1}). In other words $\phi_1^{''}$ is a solution for
the same problem as $\phi_1$ (that is \rife{011112.0}) but for the fact that
$\mu_0^2$ is replaced by $(\mu_0^2)^{''}$ in \rife{19.0.0}. Hence the argument in subsection \ref{subs:1}
applies and we obtain (see \rife{19.0.1})
$$
\int\limits_{\om_\a}\phi_1^{''}=(\mu_0^2)^{''}\int\limits_{\om_\a}\psi_0^2=
(\mu_0^2)^{''}\frac{\pi}{12\a}+\mbox{\rm O}(\a^2)=\frac{1}{6\pi\a}+\mbox{\rm O}(\a^2),
$$

and the conclusion follows by substituting this result in \rife{19.0.2}.
\fineproof

\end{document}